\documentclass[10pt]{amsart}

\usepackage{amsfonts,amsmath,latexsym,amssymb,verbatim,amsbsy,amsthm, bm, esint}
\usepackage{srcltx}
\usepackage{dsfont}

\usepackage{graphicx}

\usepackage{caption}
\usepackage{float}
\usepackage{wrapfig}

\usepackage[top=1in, bottom=1in, left=1in, right=1in]{geometry}

\usepackage[dvipsnames]{xcolor}

\usepackage{blkarray}
\usepackage{booktabs}
\usepackage[colorlinks=true, pdfstartview=FitV, linkcolor=RoyalBlue,citecolor=ForestGreen, urlcolor=blue]{hyperref}

\theoremstyle{plain}
\newtheorem{THEOREM}{Theorem}[section]

\newtheorem{theorem}[THEOREM]{Theorem}
\newtheorem{corollary}[THEOREM]{Corollary}
\newtheorem{lemma}[THEOREM]{Lemma}

\newtheorem{assumption}[THEOREM]{Assumption}
\newtheorem{notation}[THEOREM]{Notation}

\theoremstyle{definition}

\newtheorem{definition}[THEOREM]{Definition}

\theoremstyle{remark}

\newtheorem{remark}[THEOREM]{Remark}

\newtheorem{example}[THEOREM]{Example}

\newcommand{\diff}{\mathop{}\!\mathrm{d}}

\DeclareMathOperator{\diver}{div} %
\def \a {\alpha}
\def \b {\beta}
\def \g {\gamma}
\def \d {\delta}
\def \k {\kappa}
\def \e {\varepsilon}

\def \l {\lambda}
\def \n {\nabla}

\def \t {\tau}

\def \o {\omega}

\def \O {\Omega}

\def \bu {{\bf u}}
\def \bv {{\bf v}}

\def \bU {{\bf U}}

\def \cC {\mathcal{C}}

\def \cE {\mathcal{E}}

\def \cI {\mathcal{I}}

\def \cL {\mathcal{L}}
\def \cM {\mathcal{M}}

\def \cP {\mathcal{P}}
\def \cQ {\mathcal{Q}}

\newcommand{\N}{\ensuremath{\mathbb{N}}}   %%% naturals
\newcommand{\Z}{\ensuremath{\mathbb{Z}}}   %%% integers
\newcommand{\R}{\ensuremath{\mathbb{R}}}   %%% reals
\newcommand{\T}{\ensuremath{\mathbb{T}}}   %%% torus

\renewcommand{\S}{\ensuremath{\mathbb{S}}} %%% sphere

\def \lan {\langle}
\def \ran {\rangle}

\def \p {\partial}

\def \dx  {\, \mbox{d}x}
\def \dt  {\, \mbox{d}t}

\def \dy  {\, \mbox{d}y}

\def \dl   {\, \mbox{d}\lambda}

\def \dtau  {\, \mbox{d}\tau}

\def \dnu  {\, \mbox{d}\nu}

\def \ddt  {\frac{\mbox{d\,\,}}{\mbox{d}t}}

\title[WS stability and uniqueness of mv-solutions to Euler-type systems]{Weak-strong stability and uniqueness of measure-valued solutions to (nonlocal) pressureless Euler-type systems}

\author{Piotr Gwiazda}
\address{Institute of Mathematics, Polish Academy of Sciences, ul. Śniadeckich 8, 00-656 Warsaw, Poland; Interdisciplinary Centre for Mathematical and Computational Modelling, University of Warsaw, ul. Tyniecka 15/17, 02-630 Warsaw, Poland}
\email{pgwiazda@impan.pl}

\author{David N Reynolds}
\address{Faculty of Mathematics, Informatics and Mechanics of University of Warsaw, Stefana Banacha 2, 02-097 Warsaw, Poland}
\email{d.reynolds@uw.edu.pl}

\author{Agnieszka \'Swierczewska-Gwiazda}
\address{Faculty of Mathematics, Informatics and Mechanics of University of Warsaw, Stefana Banacha 2, 02-097 Warsaw, Poland}
\email{a.swierczewska-gwiaz@uw.edu.pl}

\author{Jakub Wo\'znicki}
\address{Faculty of Mathematics of University of Vienna, Oskar-Morgenstern-Platz 1, 1090 Wien, Austria }
\email{jakub.woznicki@univie.ac.at}

\begin{document}

\begin{abstract}
    We consider the pressureless Euler equation possibly coupled with non-local forces like the Poisson and alignment interactions. We overcome the lack of strict convexity in the relative entropy inequality that arises from the lack of pressure. The key to this result is the appropriate framework of the solution scheme utilizing the generalized Young measure which yields extra structural properties for the oscillation and concentration parts of the measure-valued solution. We overcome the problem of density mismatch via two complementary techniques. The first takes advantage of a Figalli--Kang-type transport inequality. The second involves establishing a generation lemma for the generalized Young measure which simultaneously establishes $L^2$-based control in flat metric between the measure-valued density and the density generating sequence. 
\end{abstract}


\maketitle

\section{Introduction}\label{S:1}

The pressureless Euler system 
    \begin{equation}\label{pES1}
\left\{
\begin{split}
\p_t\rho + \diver (\rho \bu) & = 0, \\
\p_t(\rho \bu) + \diver (\rho \bu \otimes \bu ) &=0 ,
\end{split} \right. 
\end{equation}
arises from the Zeldovich adhesion model for the formation of large-scale structure in the universe \cite{Zel}. It represents a macroscopic description of sticky particle dynamics \cite{BD2023} and can formally be derived from letting the pressure go to zero in the compressible Euler equations, or by dropping the temperature to zero in the Boltzmann equations, and more rigorously (at least locally in time) by considering monokinetic solutions to the Vlasov equation \cite{BrenierGrenier1998}.  It is further used in several numerical schemes taking advantage of splitting methods for the Euler equations \cite{ToroVazquez2012}, as well as more recent implicit-explicit (IMEX) schemes \cite{BoscheriDimarco2020}. The pressureless system \eqref{pES1} also sits as the basis for pressureless flows with nonlocal forces
\begin{equation}\label{pES-non}
\left\{
\begin{split}
\p_t\rho + \diver (\rho \bu) & = 0, \\
\p_t(\rho \bu) + \diver (\rho \bu \otimes \bu ) &=F(\rho,\bu) ,
\end{split} \right.
\end{equation}
where $$F(\rho,\bu)=-\rho\nabla \Phi_\rho, \quad -\Delta_x \Phi_\rho=\rho-M_\rho,$$ gives the pressureless Euler-Poisson system \cite{bardi2001eulerpoisson,tadmor1999critical}, and 
$$F(\rho,\bu)=\rho\int_{\O} \Psi(x,y)\rho(y)(\bu(y)-\bu(x)) \dy,$$
leads to the pressureless Euler alignment system describing macroscopic flocking dynamics \cite{CFTV10,CS2007a,CS2007b,Rbook}, or traffic flow via its identification with the Aw-Rascle model \cite{CPSZ2024,WZ2025}.

Even for smooth initial data, solutions to \eqref{pES1} and \eqref{pES-non} typically develop concentrations and discontinuities in finite time due to the crossing of characteristics \cite{BressanNguyen2014,ChaeTadmor2008,EngelbergLiuTadmor2001}. The absence of a pressure term is precisely to blame for this phenomena. Despite this fact, with certain structural assumptions, both local and global well-posedness for classical solutions to \eqref{pES1} \cite{BouchutJames1999,BrenierGrenier1998,CavallettiSedjroWestdickenberg2015,HuangWang2001,EWestdickenberg2003} and \eqref{pES-non} \cite{carrillo2009critical,CCTT2016,carrillo2014globalalignment,tadmor2011critical,tadmor1999critical,choi2026finitetimebreakdowneuleralignmentsupercritical} and even the full Euler-Poisson-Alignment system \cite{BHATNAGAR202382} can be provided.

The coexistence of strong solutions (when they exist) and highly non-unique weak solutions (due to the convex integration analysis, cf. \cite{Donatelli03072015,CFGS17}) motivates the study of weak–strong uniqueness principles. Such results assert that a suitably defined weak solution coincides with a strong solution emanating from the same initial data, as long as the latter exists. For the pressureless Euler system, this question is particularly subtle due to the possibility of oscillation and concentration effects and the lack of regularizing mechanisms from e.g. the pressure. Moreover, lack of pressure means that the standard relative entropy considerations fail in capturing the distance between the strong and weak density.

Measure-valued solutions provide a natural framework for capturing oscillations and concentrations in nonlinear PDEs, particularly in regimes where classical or even weak solutions fail to be well-posed. This approach was first systematically introduced by DiPerna \cite{DiPerna1985} and further developed with Majda \cite{DiPernaMajda1987}, who showed that weak limits of bounded sequences of approximate solutions can be described by parametrized probability measures encoding oscillatory behavior, together with concentration effects. Their framework was later refined by Alibert and Bouchitt\'e \cite{AlibertBouchitte1997} through the introduction of generalized Young measures, which incorporate an additional concentration measure and recession function to fully describe nonlinear limits. This added structure is essential in problems such as the pressureless Euler system, where solutions naturally develop concentrations (e.g., delta-shocks) and classical notions of weak convergence are insufficient. 

The weak-strong uniqueness principle was established for the incompressible Euler equations by Brenier et al. \cite{BrenierDeLellisSzekelyhidi2011}, using the relative entropy (or modulated energy) method introduced by Dafermos \cite{Dafermos1979} for conservation laws. It was subsequently extended to compressible isentropic Euler and general hyperbolic conservation laws by Gwiazda et al. \cite{AIHPC_2020__37_3_683_0,GwiazdaSwierczewskaWiedemann2015}, to polyconvex elastodynamics by Demoulini et al. \cite{Demoulini2012}, and to Navier-Stokes by Feireisl et al. \cite{FeireislGwiazdaSwierczewskaWiedemann2016b}.

A key conceptual advance, advocated in the framework of \cite{BrenierDeLellisSzekelyhidi2011} and subsequently adopted in the compressible setting \cite{AIHPC_2020__37_3_683_0,GwiazdaSwierczewskaWiedemann2015}, is to define measure-valued solutions intrinsically, as objects independent of any generating sequence, rather than as limits of particular approximations. This approach is essential for the logical structure of weak-strong uniqueness. If the solution concept depended on the choice of approximating sequence, then uniqueness would be a statement about a specific approximation scheme rather than a property of the solution itself. Indeed, in the compressible Euler case it has been shown that not every measure-valued solution can be approximated by a sequence of weak solutions \cite{CFKW2015}. Further, we provide Example \ref{example} of a measure-valued solution to \eqref{pES1} which cannot be generated by an approximating sequence of bounded weak solutions, further justifying our motivation for obtaining weak-strong uniqueness and stability properties which do not rely on approximating sequences of weak solutions. 

The pressureless Euler system presents substantial obstacles that prevent a direct application of the standard relative entropy approach. The relative entropy functional for isentropic Euler is modulated by the pressure potential $P(\rho) = \rho^\gamma / (\gamma-1)$, whose strict convexity is the key mechanism for absorbing the density discrepancy between the measure-valued and strong solutions. In the pressureless case ($P \equiv 0$), this mechanism breaks down entirely, and the density error cannot be controlled by energetic means alone. As observed in \cite{WZ2025}, previous works on measure-valued weak-strong uniqueness required either Lebesgue or dual Sobolev norm control on the density, which is not readily available in these nonlocal pressureless systems.

A further complication is the ill-posedness of the pressureless system at the level of bounded admissible weak solutions in dimensions $d \geq 2$. Within \cite{CFGS17} the authors constructed, via convex integration, infinitely many global dissipative weak solutions for bounded initial data. More recently, non-uniqueness of admissible weak solutions to the two-dimensional Riemann problem has been demonstrated \cite{HSW25}. These results seem to imply that rather than the class of admissible weak solutions, the measure-valued framework with an intrinsic admissibility criterion may be a more appropriate setting in which to formulate the weak-strong uniqueness theory.

A crucial ingredient that makes the pressureless problem tractable is the Wasserstein transport inequality of Figalli and Kang \cite{FK19}, introduced in their rigorous derivation of the pressureless Euler system with nonlocal alignment as the hydrodynamic limit of the kinetic Cucker-Smale model \cite{CS2007a}. The central difficulty in that analysis is that the entropy of the pressureless limit system lacks strict convexity in the density variable, which they overcome by combining relative entropy quantities with the 2-Wasserstein distance. This approach has since been recognized as a key tool in the pressureless setting; see also \cite{CarrilloChoiJinwook2021,ChoiSalem2017,HaKimZhang2018,KangHaKimShim2020,WZ2025}.

In this work we prove a general weak-strong uniqueness theorem for generalized measure-valued solutions to the multi-dimensional pressureless Euler equations, with or without the nonlocal forces, as in \eqref{pES-non}. Our solution concept is formulated using a generalized Young measure which takes into account oscillations and concentrations via a recession function. The added structure of the generalized Young measure allows us to resolve a gap in the existing theory, in particular when applied to the Euler-Poisson equations as in \cite{CDGS2024}, we avoid any implicit dependence on the generating sequence, providing a weak-strong uniqueness principal independent of the approximation sequence.

Our main technical contribution is twofold; both regarding the need for extra control over the distance between measure-valued and strong densities that arises from the lack of density control in the standard relative entropy argument. The first takes advantage of a Figalli-Kang type Wasserstein-inequality to bound the $W_2-$distance by the relative entropy. The second involves a parallel construction of a generating sequence for the generalized Young measure that simultaneously tracks the distance between the density generating functions and the measure-valued density in $L^2(0,T;(\cM^+(\T^d),d_f))$, which is necessary for controlling the nonlocal alignment forces which involve products of densities and velocities across different spatial coordinates. In particular, this second tool, which we place in Appendix \ref{App1} for ease of reading, is crucial as the only aspect of the PDE used is the conservation of mass of the measure-valued solution, and is not a sequence of approximating weak solutions.

We further provide a weak-strong stability result in terms of the discrepancy between the densities and momenta in total variation. In particular, this provides that even in the case where initial datum leads to finite time blowup of the corresponding smooth solution, that the measure-valued solution remains close to its ``nearest'' smooth solutions which do retain global regularity.

The paper is outlined as follows. Section \ref{S:2} gives the relevant preliminaries and establishes the notation for generalized Young measures and the statement of the main results (see Theorem \ref{t:main}). Section \ref{S:3} establishes the necessary Figalli-Kang type inequality for controlling density mismatch. Section \ref{S:4} establishes the relative entropy inequality which provides the weak-strong stability result. Section \ref{S:5} ultimately establishes the weak-strong uniqueness result. A necessary lemma for closing the relative entropy inequality in the presence of alignment forces is provided in Appendix \ref{App1}.

\section{Preliminaries and statement of main results}\label{S:2}
We consider the pressureless Euler system with nonlocal forces
    \begin{equation}\label{pES}
\left\{
\begin{split}
\p_t\rho + \diver (\rho \bu) & = 0, \\
\p_t(\rho \bu) + \diver (\rho \bu \otimes \bu ) &= -a \rho \nabla_x \Phi_\rho+b\rho\int_{\O} \Psi(x,y)\rho(y)(\bu(y)-\bu(x))\dy,\\
-\Delta_x\Phi_\rho&=\rho-M_\rho,
\end{split} \right. 
\end{equation}
where $M_\rho=0$ if $\O=\R^d$ and $M_\rho=\frac{1}{|\T^d|}\int_{\T}\rho \dx$ if $\O=\T^d$. The constants $a,b\geq0$ represent whether or not the Euler-Poisson and/or/nor Euler-Alignment nonlocal forces are active. In particular, $a=b=0$ returns pressureless Euler, while $a>0, b=0$ gives Euler-Poisson, $a=0, b>0$ the Euler-Alignment system, and $a,b>0$ a combined Euler-Poisson-Alignment system. The techniques developed here apply to all of these systems and hence we provide a unified theory controlling all of the terms at once.

One important consideration is that of the domain $\Omega$. In principal, the techniques developed here are applicable to smooth bounded domains with the impermeability condition and Neumann boundary conditions
\begin{align*}
    \bu\cdot n=0, \quad \n_x\Phi_\rho\cdot n=0, \quad \text{on} \quad \p\O\times (0,T),
\end{align*}
however, questions as to the persistence of strong solutions through interactions with the boundary, and the physical interpretation of a pressureless macroscopic system on a bounded domain leads us to focus at present on the periodic domain $\O=\T^d$, and the whole space $\O=\R^d$. Further, for technical reasons (see Appendix \ref{App1} for the lemma about the generation of generalized Young measures), if the alignment force is active $b>0$, we restrict ourselves further to only the periodic domain $\O=\T^d$.

Further, in order to formulate this system in the measure-valued framework, it is necessary to rewrite the Poisson interaction in divergence form. Given
\begin{align*}
-\Delta_x \Phi_\rho = \rho-M_\rho,
\end{align*}
 the Poisson interaction can be expressed as
\begin{align*}
(\rho-M_\rho) \nabla_x \Phi_\rho
=
\diver
\left(
-\nabla_x \Phi_\rho \otimes \nabla_x \Phi_\rho
+ \frac{1}{2} |\nabla_x \Phi_\rho|^2 I
\right).
\end{align*}
This reformulation avoids ill-defined products of weakly converging quantities and allows the nonlocal interaction to be treated within the framework of generalized Young measures. Although also a nonlocal force, the alignment operator
\begin{align*}
\rho\mathcal A[\rho,\bu](x)
=
\rho\int_{\O}
\Psi(x,y)\rho(y)(\bu(y)-\bu(x))\dy
\end{align*}
is bilinear, and thus does not require any further reformulation. For the purposes of this work, we consider the following assumptions for interaction kernel $\Psi(x,y)$
\begin{assumption}\label{ass:on_kernel_psi}
    We assume that $\Psi \in W^{1,\infty}$ and it is an even, continuous, function with values in positive semidefinite matrices.
\end{assumption}
Thus we rewrite the whole system as
\begin{equation}\label{pESdiv}
\left\{
\begin{split}
\p_t\rho + \diver (\rho \bu) & = 0, \\
\p_t(\rho \bu) + \diver (\rho \bu \otimes \bu ) &= -\frac{a}{2} \nabla_x|\n_x\Phi_\rho|^2+a\diver[ \n_x\Phi_\rho \otimes \n_x\Phi_\rho]-aM_\rho\n_x\Phi_\rho\\
& \ \ \ \ \ +b\rho\int_{\O} \Psi(x,y)\rho(y)(\bu(y)-\bu(x))\dy, \qquad a,b\geq0,\\
-\Delta_x\Phi_\rho&=\rho-M_\rho.
\end{split} \right. 
\end{equation}

The natural energy associated with the system is given by
\begin{align*}
\cE(t)
=
\int_{\O}
\frac{1}{2} \rho |\bu|^2\dx + \frac{a}{2}
\int_{\O}
|\nabla_x \Phi|^2 \dx,
\end{align*}
where the first term represents the kinetic energy and the second corresponds to the interaction energy induced by the Poisson potential. Formally differentiating in time yields
\begin{align*}
\ddt \cE(t)
=
-\frac{b}{2}
\iint_{\O \times \O}
\Psi(x,y)\rho(x)\rho(y)|\bu(x)-\bu(y)|^2\dx\dy \leq 0.
\end{align*}
Thus, the Poisson interaction is conservative, while the alignment mechanism, when $b>0$, induces dissipation.

The relative entropy between two solutions $(\rho,\bu)$ and $(r,\bU)$ is thus defined
\begin{align*}
     \cE_{rel}(\t)=\int_{\O}\frac{1}{2}\rho|\bu-\bU|^2 \dx + \int_\Omega \frac{a}{2}|\nabla_x\Phi_{\rho} - \nabla_x\Phi_r|^2\diff x,
\end{align*}
where we can see, especially in the case $a=0$, that vanishing of the relative entropy fails to provide information about the densities $\rho$ and $r$. This is the main hurdle to overcome, however, we first must introduce the relevant notation for defining our notion of solution to the system \eqref{pESdiv}.

\subsection{Generalized Young measures}
In here, we shall recall the notion of the generalized Young measures, which will be of help in defining a proper notion of a solution. First introduced by DiPerna and Majda in \cite{DiPerna1985,DiPernaMajda1987} and later refined by Alibert and Bouchitt\'e \cite{AlibertBouchitte1997}, this tool is of great importance when dealing with equations with very weak a priori estimates such as pressureless Euler equations, in which one can hope only for $L^\infty(0, T; L^1(\Omega))$ bounds for density and momentum. To fix our attention for a moment let $Y\subset\R^N$ be a measurable subset and $y$ be its elements. Moreover, let $\{(a_k(y), b_k(y))\}_{k\in\N}$ be a sequence of functions such that $\{a_k\}_{k\in\N}$ is bounded in $L^2(Y,\R^l)$ and $\{b_k\}_{k\in\N}$ in $L^2(Y, \R^m)$. Then, define a sphere
$$
\S^{l+m-1}_2 := \{(a_\infty, b_\infty)\in \R^{l+m}:\, |a_\infty|^4 + |b_\infty|^4 = 1\}.
$$
Having the above, one may follow \cite{AlibertBouchitte1997} to deduce, that there exists a triple
\begin{itemize}
    \item $\nu\in L^\infty_w(Y; \cP(\R^{l+m}))$,
    \item $\lambda\in \cM^+(\overline{Y})$,
    \item $\nu^\infty\in L^\infty_w((Y, \lambda); \cP(\S^{l+m-1}_2))$,
\end{itemize}
such that for any Carath\'{e}odory $f$ for which the recession function
$$
f^\infty(y, a_\infty, b_\infty) := \lim_{\substack{y'\to y, \\ (a_\infty', b_\infty')\to (a_\infty, b_\infty),\\ r\to +\infty}}\frac{f(y', r^2a_\infty', r^2b_\infty')}{r^4}
$$
is well-defined and continuous on $\overline{Y}\times \S^{l+m-1}_2$ the following limit holds
\begin{align*}
    f(y, a_k(y), b_k(y))\dy \overset{*}{\rightharpoonup}\int_{\R^{l+m}}f(y, a, b)\dnu_y(a, b)\dy + \int_{\S^{l+m-1}_2}f^\infty(y, a_\infty, b_\infty)\dnu^\infty_y(a_\infty, b_\infty)\diff\lambda(y),
\end{align*}
where the convergence is understood in the sense of weak* in measures. To deal with convolutions, one can easily verify that generalized Young measures tensorize. That is, if we have two sequences $\{a_k(y)\}_{k\in\N}$ and $\{a_k(z)\}_{k\in\N}$, then with the established notation above, the sequence $\{(a_k(y), a_k(z))\}_{k\in\N}$ generates a Young measure $\nu_{y, z}^2 := \nu_{y}\otimes \nu_z$ in terms of oscillation, and $\nu_{y,z}^{2,\infty} := \nu_{y}^\infty\otimes \nu_z^\infty$, $\lambda^2(\diff y, \diff z) := \lambda(\diff y)\otimes \lambda(\diff z)$ in terms of concentration. Those products have following properties
\begin{enumerate}
    \item \textbf{Symmetry}: 
    \begin{multline*}
        \int_{\R^{2l}}f(a, a')\dnu^2_{y,z}(a, a')\dy\diff z + \int_{(\S^{l-1}_2)^2}f^\infty(a_\infty, a'_\infty)\dnu^{2,\infty}_{y,z}(a_\infty, a'_\infty)\diff\lambda^2(y, z)\\ = \int_{\R^{2l}}f(a', a)\dnu^2_{z,y}(a', a)\dy\diff z + \int_{(\S^{l-1}_2)^2}f^\infty(a'_\infty, a_\infty)\dnu^{2,\infty}_{z,y}(a'_\infty, a_\infty)\diff\lambda^2(z, y),
    \end{multline*}
    for all admissible functions $f$ (that is continuous, vanishing at infinity, whose recession functions are well defined).
    \item \textbf{Consistency}: If $f := f(a, a')$ is a function of one variable, that is there exists $g$, such that $f(a, a') = g(a)$, then
    \begin{multline*}
        \int_{\R^{2l}}f(a, a')\dnu^2_{y,z}(a, a')\dy\diff z + \int_{(\S^{l-1}_2)^2}f^\infty(a_\infty, a'_\infty)\dnu^\infty_{y,z}(a_\infty, a'_\infty)\diff\lambda^2(y, z)\\ = \int_{\R^{l}}g(a)\dnu_{y}(a)\dy + \int_{\S^{l-1}_2}g^\infty( a_\infty)\dnu^\infty_{y}(a_\infty)\diff\lambda(y),
    \end{multline*}
    for all admissible $f$ and $g$.
\end{enumerate}
Although these sequences are not the major focus in the present work, this setting presents a motivation and a clear approximation path which could be established to show existence of measure-valued solutions to general systems, as in \cite{CDGS2024} which is particularly relevant for the case at hand. We choose to detach ourselves due to the work of \cite{CFKW2015} which suggests that there exist measure-valued solutions which are not obtained from sequences of weak solutions. We would also note that what we have just presented is quite general and will heavily simplify when we come back to the example of pressureless Euler equations below. As a last, but not least, point we would like to note that in our case of a time-dependent equation, $Y$ will be of the form $(0, T)\times \Omega$. In this case we can find a decomposition
$$
\lambda = \lambda_t(\diff x)\otimes \diff t,
$$
where $t\mapsto \lambda_t$ is a bounded and measurable function going from $(0, T)$ to $\cM^+(\overline\Omega)$ cf. \cite{BrenierDeLellisSzekelyhidi2011}.

\subsection{Definition of a solutions and the main result}

To begin our simplification, let us introduce the variables that we are going to be working with, which may be seen as emanating from $(\sqrt{\rho}, \sqrt{\rho}\bu, \nabla_x\Phi_\rho)$. Since any approximating sequence of this form should be in principle bounded in $L^\infty(0, T; L^2(\Omega))$ it is reasonable to assume we have an appropriate triple $(\nu_{t, x}, \lambda, \nu^\infty_{t,x})$ such that
\begin{itemize}
    \item $\nu\in L^\infty_w((0, T)\times \Omega; \cP(\R_+\times\R^d\times \R^d))$,
    \item $\lambda\in \cM^+([0, T]\times \overline{\Omega})$,
    \item $\nu^\infty\in L^\infty_w(((0, T)\times \Omega, \lambda); \cP(\S^{2d}_2))$.
\end{itemize}
Throughout the article, we shall use the following dummy variables
\begin{align}
    &(s,\b, F) \sim (\sqrt{\rho},\sqrt{\rho}\bu, \nabla_x\Phi_\rho) \quad \text{in} \ \nu_{t,x},\label{def:dummy_nu}\\
    &(s_\infty,\b_\infty, F_\infty) \sim (\sqrt{\rho},\sqrt{\rho}\bu, \nabla_x\Phi_\rho) \quad \text{in} \ \nu_{t,x}^\infty,\nonumber
\end{align}
as well as whenever we deal with convolution and need to indicate separate variables in space
\begin{align*}
    &(s',\b') \sim (\sqrt{\rho},\sqrt{\rho}\bu) \quad \text{in} \ \nu_{t,y},\\
    &(s'_\infty,\b'_\infty) \sim (\sqrt{\rho},\sqrt{\rho}\bu) \quad \text{in} \ \nu_{t,y}^\infty,
\end{align*}
and it will become clear shortly what that actually means. To further simplify, we define a standard bracket notation for an action of a measure on a function, that is
$$
\int_{\R^{2d+1}}f(t, x, s, \beta, F)\dnu_{t, x}(s, \beta, F)\diff x\diff t = \langle \nu_{t,x}, f(t, x, s, \beta, F)\rangle\diff x\diff t,
$$
and
$$
\int_{\S^{2d}_2}f^\infty(t, x, s_\infty, \beta_\infty, F_\infty)\dnu^\infty_{t, x}(s_\infty, \beta_\infty, F_\infty)\diff\lambda_t(x)\diff t = \langle \nu^\infty_{t,x}, f^\infty(t, x, s_\infty, \beta_\infty, F_\infty)\rangle \diff\lambda_t(x)\diff t,
$$
as well as the less standard notation
$$
\overline{f}\diff x\diff t := \langle \nu_{t,x}, f(t, x, s, \beta, F) \rangle\diff x\diff t + \langle \nu^\infty_{t,x}, f^\infty(t, x, s_\infty, \beta_\infty, F_\infty)\rangle \diff\lambda_t(x)\diff t.
$$
Similarly, we define those quantities for a tensorized Young measures $\nu^2_{t, x, y}$ and $\nu^{2,\infty}_{t,x,y}$ whenever we deal with convolution. In this spirit, we can define our variables, which will appear in the notion of the solution.
\begin{notation}\label{not:bar_variables}
We define the following quantities.
    \begin{align*}
    &\overline{\rho}=\langle \nu_{t,x},s^2\rangle+\langle \nu_{t,x}^\infty,s^2_\infty\rangle\l,\\
    &\overline{\rho\bu}=\langle \nu_{\tau,x},s\,\beta\rangle+\langle \nu_{\tau,x}^\infty,s_\infty\beta_\infty\rangle\l,\\
    &\overline{\rho\bu\otimes \bu}=\langle \nu_{\tau,x},\beta \otimes \beta\rangle+\langle \nu_{\tau,x}^\infty,\beta_\infty\otimes \b_\infty\rangle\l,\\
    &\overline{\rho|\bu|^2}=\langle \nu_{\tau,x},|\b|^2\rangle+\langle \nu_{\tau,x}^\infty,|\b_\infty|^2\rangle\l,\\
    &\overline{\nabla_x\Phi_\rho}=\langle \nu_{t,x},F\rangle,\\
    &\overline{|\nabla_x\Phi_\rho|^2}=\langle \nu_{t,x},|F|^2\rangle+\langle \nu_{t,x}^\infty,|F_\infty|^2\rangle\l,\\
    &\overline{\nabla_x\Phi_\rho\otimes\nabla_x\Phi_\rho}=\langle \nu_{t,x},F\otimes F\rangle+\langle \nu_{t,x}^\infty,F_\infty\otimes F_\infty\rangle\l.\\
    &\overline{\rho(x)\rho(y)(\bu(x) - \bu(y))} = \langle \nu^2_{t,x,y},ss'(s'\beta - s\beta') \rangle+\langle \nu_{t,x,y}^{2,\infty},s_\infty s'_\infty(s'_\infty \beta_\infty - s_\infty\beta'_\infty)\rangle\l^2\\
    &\overline{\rho(x)\rho(y)(\bu(x) - \bu(y))\Psi(x-y)(\bu(x) - \bu(y))} = \langle\nu^2_{t,x, y}, (s'\beta - s\beta')\Psi(x-y)(s'\beta - s\beta')\rangle\\
    &\qquad\qquad\qquad\qquad\qquad\qquad\qquad\qquad\qquad\qquad+ \langle\nu^{2,\infty}_{t,x, y}, (s_\infty'\beta_\infty - s_\infty\beta'_\infty)\Psi(x-y)(s'_\infty\beta_\infty - s_\infty\beta'\infty)\rangle\l^2
\end{align*}
\end{notation}

For further clarification, let us show how to quickly unravel this notation for the convenience of the reader. Take for example $\rho\bu$. Since we work with the variables $(\sqrt{\rho}, \sqrt{\rho}\bu)$ we need to find them in this term; we can easily write
$$
\rho\bu = \sqrt{\rho}\sqrt{\rho}\bu.
$$
Then, we identify $s$ as a dummy variable for $\sqrt{\rho}$ and $\beta$ as a dummy variable for $\sqrt{\rho}\bu$ cf. \eqref{def:dummy_nu}. Hence, under the measure $\nu_{t,x}$ we need to find the term $s\beta$. Similarly, under $\nu_{t, x}^\infty$ we should have $s_\infty\beta_\infty$. Thus,
$$
\overline{\rho\bu}=\langle \nu_{\tau,x},s\,\beta\rangle+\langle \nu_{\tau,x}^\infty,s_\infty\beta_\infty\rangle\l.
$$
We are now ready to define the measure-valued solution to the equations \eqref{pESdiv}.

\begin{definition}\label{def:measure_solutions}
    We will say that the triple $(\nu, \lambda, \nu^\infty)$ is a measure-valued solution to the system \eqref{pESdiv} with measure-valued initial datum $(\nu_0, \lambda_0, \nu^\infty_0)$ iff
    \begin{itemize}
        \item $\nu\in L^\infty_w((0, T)\times \Omega; \cP(\R^+\times \R^{d}\times \R^d))$,
        \item $\lambda\in \cM^+([0, T]\times \overline{\Omega})$,
        \item $\nu^\infty\in L^\infty_w(((0, T)\times \Omega, \lambda); \cP(\S^{+}_2))$,
    \end{itemize}
    where  
    $$
    \S^+_2 := \{(s_\infty, \beta_\infty, F_\infty)\in \R^{2d+1},\, |s_\infty|^4 + |\beta_\infty|^4 + |F_\infty|^4= 1,\, s_\infty\geq 0\},
    $$
    and if $b=0$, $\O=\T^d$ or $\O=\R^d$ and if $b>0$, then $\O=\T^d$, to distinguish we let $c(\O)=1$ when $\O=\R^d$ and $c(\O)=0$ if $\O=\T^d$.
    
    Further, the following relations are satisfied
    \begin{itemize}
    \item \textit{Continuity equation:}
    \begin{equation}\label{CE-mv}
    \begin{split}
        \int_\O \overline{\rho}(\tau, x)\psi(\tau,x) \dx- \int_\O \overline{\rho}_0\psi(0,x) \dx =\int_0^\t\int_\O \overline{\rho}(t, x) \p_t\psi+ \overline{\rho\bu}(t, x) \cdot \n_x\psi \dx \dt,
    \end{split}
    \end{equation}
    for a.e. $\t\in(0,T)$ and every $\psi \in C^1([0,T)\times \overline{\O})$.
    \item \textit{Momentum equation:}
    \begin{equation}\label{ME-mv}
    \begin{split}
        \int_\O &\overline{\rho\bu}(\tau, x) \cdot \phi(\t,x) \dx-\int_\O \overline{\rho\bu}_0 \cdot \phi(0,x) \dx \\
        &=\int_0^\t\int_\O \overline{\rho\bu}(t, x)\cdot \p_t\phi \dx\dt+\int_0^\t\int_\O \overline{\rho\bu\otimes\bu}: \n_x \phi \dx\dt  + \int_{0}^{\tau}\int_{\Omega}^{}\frac{a}{2}\overline{|\nabla_x\Phi_{\overline{\rho}}|^2}\operatorname{div}_x\phi\diff x\diff t\\
        &\quad \quad+ \int_0^\tau\int_{\Omega}a\overline{\nabla_x\Phi_{\overline{\rho}}\otimes\nabla_x\Phi_{\overline{\rho}}}:\nabla_x\phi\diff x\diff t - \int_0^\tau\int_{\Omega}a M_{\overline{\rho}}\overline{\nabla_x\Phi_{\overline{\rho}}}\cdot\phi\,\diff x\diff t\\
        &\quad\quad -\frac{b}{2}\int_0^\tau\int_{\Omega\times\Omega}\overline{\rho(x)\rho(y)(\bu(x) - \bu(y))}\Psi(x-y)(\phi(x) - \phi(y))\diff y\diff x\diff t
    \end{split}
    \end{equation}
    for a.e. $\t \in (0,T)$ and every $\phi \in C^1([0,T)\times \overline{\O}; \R^d)$.
    \item \textit{Poisson equation:} For a.e.~$\tau\in(0,T)$ and all $\chi\in C^1(\overline{\Omega})$ it holds that
	\begin{align}
	    \int_{\Omega}^{}\langle\nu_{(\tau,\cdot)},F\rangle\cdot\nabla_x\chi\,\diff x=\int_{\Omega}^{}\overline{\rho}\,\chi\,\diff x-M_{\overline{\rho}}\,\int_{\Omega}^{}\chi\,\diff x.\label{eq:poissonmvs}
	\end{align}
    Moreover, there exists some $\Phi_{\overline{\rho}}\in L^2\big((0,T);W^{1,2}\cap W^{1,q}(\Omega)\big)$ for $1 < q < \frac{d}{d-1}$, such that
	\begin{align}
	    \nabla_x\Phi_{\overline{\rho}}=\langle\nu,F\rangle\label{eq:nablaPhiidentity}
	\end{align}
	a.e.~on $(0,T)\times \Omega$.
    \item \textit{Energy inequality:}
    \begin{equation}\label{EN-mv}
    \begin{split}
        &\int_\O \overline{\rho|\bu|^2}(\tau, x)\dx + \int_{\Omega}^{}a\overline{|\nabla_x\Phi_{\overline{\rho}}|^2}(\tau, x)\diff x +\frac{b}{2}\int_0^\tau\int_{\Omega\times \Omega}\overline{\rho(x)\rho(y)(\bu(x) - \bu(y))\Psi(x-y)(\bu(x) - \bu(y))}\diff y\diff x\diff t\\
        &\qquad\leq \int_\O \overline{\rho|\bu|^2}_0\dx + \int_{\Omega}^{}a\overline{|\nabla_x\Phi_{\overline{\rho}}|^2}_0\diff x
    \end{split}
    \end{equation}
    for a.e. $\t \in (0,T)$.
    \end{itemize}
\end{definition}
For the convenience of the reader, let us unravel one of the relations stated in the definition above to its full form. For example, for the continuity equation \eqref{CE-mv} using the Notation \ref{not:bar_variables} we can write equivalently
\begin{equation*}
    \begin{split}
        &\int_\O \langle\nu_{\tau, x}, s^2\rangle\psi(\tau,x) \dx + \int_\Omega \langle\nu_{\tau, x}^\infty, s_\infty^2\rangle\psi(\tau, x)\diff\lambda_\tau(x) - \int_\O \langle\nu_{0, x}, s^2\rangle\psi(0,x) \dx - \int_\Omega \langle\nu_{0, x}^\infty, s_\infty^2\rangle\psi(0, x)\diff\lambda_0(x) \\
        &=\int_0^\t\int_\O \langle\nu_{t, x}, s^2\rangle\p_t\psi+ \langle\nu_{t, x}, s\beta\rangle\cdot \n_x\psi \dx \dt + \int_0^\tau\int_\Omega\langle\nu_{t, x}^\infty, s_\infty^2\rangle\p_t\psi + \langle\nu_{t, x}^\infty, s_\infty\beta_\infty\rangle\cdot\nabla_x\psi\diff\lambda_t(x)\diff t.
    \end{split}
\end{equation*}

We provide here a quick example to further motivate considering measure-valued solutions, independent of an approximating solution scheme.

\begin{example}\label{example}
    Consider the homogeneous Young measure $$\nu_{t,x}(s, \beta)=\d_{\{1\}}(s) \otimes \left(\frac{1}{2}\d_{\{-1\}}+\frac{1}{2}\d_{\{1\}}\right)(\beta).$$ Then $\nu_{t,x}$ solves \eqref{pES1}, but cannot be generated by a sequence of bounded weak solutions.
    \begin{proof}
        It is easy to see that $\nu_{t,x}$ is a measure solution, as it is homogeneous (most of the integrals will be equal simply to $0$). To see that there is no possibility of generation, by way of contradiction assume that there exists a generating sequence $(\rho_n, u_n)$ of bounded weak solutions to \eqref{pES1}. Since the Young measure on the dummy variable $s\sim \sqrt{\rho}$ is equal to $\delta_{\{1\}}$, by standard properties of Young measures we have
        $$
        \sqrt{\rho_n}\rightarrow 1,\text{ in measure,} 
        $$
        hence also
        $$
        \rho_n \rightarrow 1, \text{ in measure.}
        $$
        Moreover, it is easy to check that the Young measure generated by $\rho_n u_n^2$ is also equal to $\delta_{\{1\}}$, which again implies
        $$
        \rho_n u_n^2\rightarrow 1,\text{ in measure.}
        $$
        Furthermore, as the sequences are bounded, by Vitali's theorem 
        $$
        \rho_n, \rho_nu_n^2\rightarrow 1,\text{ strongly in }L^p,\quad 1\leq p < +\infty.
        $$
        Thus, using the continuity equation
        $$
        \partial_x(\rho_n u_n) = -\partial_t(\rho_n - 1)\rightarrow 0\text{, strongly in } W^{-1, p}((0, T)\times \T^d),\quad 1\leq p < +\infty.
        $$
        and from the momentum equation
        $$
        \partial_t(\rho_n u_n) = -\partial_x(\rho_n u_n^2)\rightarrow 0\text{, strongly in }W^{-1, p}((0, T)\times \T^d),\quad 1\leq p < +\infty.
        $$
        Hence, by the negative-norm Poincar\'{e}'s inequality
        \begin{align}\label{ineq:negative_poincare_example}
        \left\|\rho_n u_n - \frac{1}{T|\T^d|}\int_{0}^T\int_{\T^d}\rho_n u_n\diff x\diff t\right\|_{L^p((0, T)\times\T^d)}\leq C\|\nabla_{t,x}(\rho_n u_n)\|_{W^{-1, p}((0, T)\times\T^d)} \rightarrow 0.
        \end{align}
        At the same time from the form of the Young measure
        \begin{align*}
            \rho_n u_n \rightharpoonup \int_{\R^2}s\beta\diff\delta_{\{1\}}(s)\diff\left(\frac{1}{2}\d_{\{-1\}}+\frac{1}{2}\d_{\{1\}}\right)(\beta) = 0,\text{ weakly in }L^1((0, T)\times\T^d)
        \end{align*}
        thus
        $$
        \frac{1}{T|\T^d|}\int_{0}^T\int_{\T^d}\rho_n u_n\diff x\diff t \rightarrow 0,\text{ strongly in }L^p((0, T)\times \T^d),\quad 1\leq p < +\infty,
        $$
        and in turn by \eqref{ineq:negative_poincare_example}
        \begin{align*}
            \rho_n u_n\rightarrow 0,\text{ strongly in }L^p((0, T)\times \T^d),\quad 1\leq p < +\infty.
        \end{align*}
        But, $\rho_n$ and $\rho_n u_n^2$ already converge strongly in $L^p((0, T)\times \T^d),\,1\leq p < +\infty$, by Vitali's theorem to $1$, hence their product also converges to $1$, as they are bounded, while strong convergence of $\rho_n u_n$ to $0$ implies the strong convergence of $\rho_n^2 u_n^2$ to $0$, which is a contradiction.
    \end{proof}
\end{example}

 As noted, generating sequences of the generalized Young measures is not the focus of this work, but we will need the following lemma in order to control the relative entropy inequality in the case when the nonlocal alignment force is present.

\begin{lemma}\label{lem:alignment_inequality_for_relative_entropy}
    Suppose $(\nu, \lambda, \nu^\infty)$ is a measure-valued solution given by the Definition \ref{def:measure_solutions} with spatial domain being $\T^d$, and $(r, \bU)\in L^1((0, T)\times \Omega)\times L^1(0, T; W^{1,\infty}(\Omega))$ are given functions, as well as that $\Psi$ satisfies Assumption \ref{ass:on_kernel_psi}. Then, there exists a constant $C > 0$ such that for a.e. $\tau\in (0, T)$
    \begin{equation*}
        \begin{split}
            &\Bigg|\int_0^\tau\int_{\T^d\times\T^d}\langle\nu^2_{t,x, y}, (s^2\bU(x) - s\beta)((s')^2 - r(y))\rangle\Psi(x-y)(\bU(x) - \bU(y))\diff y\diff x\diff t\\
        &\qquad + \int_0^\tau\int_{\T^d\times\T^d}\langle\nu^{2,\infty}_{t,x, y}, (s^2_\infty\bU(x) - s_\infty\beta_\infty)((s'_\infty)^2 - r(y))\rangle\Psi(x-y)(\bU(x) - \bU(y))\diff\l(y)\diff\l(x)\diff t\Bigg|\\
        &\leq C\|\Psi\|_{W^{1,\infty}}\|\bU\|_{L^\infty_t W^{1,\infty}_x}\left(\int_0^\tau\int_{\T^d}\overline{\rho|\bu - \bU|}^2\diff x\diff t + \int_0^\tau d_f^2(\overline{\rho}_t, r_t)\diff t\right).
        \end{split}
    \end{equation*}
\end{lemma}

\begin{proof}
    Let $\vartheta_n = (s_n, \beta_n, F_n)$ be a sequence of functions generating $(\nu, \lambda, \nu^\infty)$ given by Lemma \ref{Appendix1}. Then, using Young's inequality
    \begin{equation}\label{ineq:generation_lemma}
        \begin{split}
        &\Bigg|\int_0^\tau\int_{\T^d\times\T^d}(s^2_n(x)\bU(x) - s_n(x)\beta_n(x))(s_n^2(y) - r(y))\Psi(x-y)(\bU(x) - \bU(y))\diff y\diff x\diff t\Bigg|\\
        &\leq C\|\Psi\|_{W^{1,\infty}}\|\bU\|_{W^{1,\infty}}\int_0^\tau\int_{\T^d}|s^2_n(x)\bU(x) - s_n(x)\beta_n(x)|d_f(s_n^2\diff y, r\diff y)\diff t\diff x\\
        &\leq C\|\Psi\|_{W^{1,\infty}}\|\bU\|_{W^{1,\infty}}\left(\int_0^\tau\int_{\T^d}|\beta_n - s_n\bU|^2\diff x\diff t + \int_0^\tau d^2_f(s_n^2\diff y, r\diff y)\diff t\right)
        \end{split}
    \end{equation}
    Since $\vartheta_n$ generates our measure we have
    \begin{equation*}
        \begin{split}
            &\int_0^\tau\int_{\T^d\times\T^d}(s^2_n(x)\bU(x) - s_n(x)\beta_n(x))(s_n^2(y) - r(y))\Psi(x-y)(\bU(x) - \bU(y))\diff y\diff x\diff t\\
            &\rightarrow \int_0^\tau\int_{\T^d\times\T^d}\langle\nu^2_{t,x, y}, (s^2\bU(x) - s\beta)((s')^2 - r(y))\rangle\Psi(x-y)(\bU(x) - \bU(y))\diff y\diff x\diff t\\
            &\qquad\qquad +\int_0^\tau\int_{\T^d\times\T^d}\langle\nu^{2,\infty}_{t,x, y}, (s^2_\infty\bU(x) - s_\infty\beta_\infty)((s'_\infty)^2 - r(y))\rangle\Psi(x-y)(\bU(x) - \bU(y))\diff\l(y)\diff\l(x)\diff t,
        \end{split}
    \end{equation*}
    as well as
    \begin{equation*}
        \begin{split}
            \int_0^\tau\int_{\T^d}|\beta_n - s_n\bU|^2\diff x\diff t\rightarrow \int_0^\tau\int_{\T^d}\overline{\rho|\bu - \bU|}^2\diff x\diff t.
        \end{split}
    \end{equation*}
    Moreover, by \eqref{conv:L_2_step_function_generation}
    $$
    \int_0^\tau d^2_f(s_n^2\diff y, r\diff y)\diff t \rightarrow \int_0^\tau d_f(\overline{\rho}_t, r\diff y)\diff t.
    $$
    Hence, converging in both sides of \eqref{ineq:generation_lemma} gives our thesis.
\end{proof}

Of particular importance is that the generation lemma of Appendix \ref{App1} does not use a sequence of weak solutions to an approximating scheme, but rather arbitrary functions which provides an extra property which allows for further control of the limiting measure-valued density.\\

We now define the measure-valued energy and relative entropy as usual
\begin{align}
    \cE^{mv}(\tau)=\int_\O\frac{1}{2}\overline{\rho|\bu|^2} \dx + \int_\Omega \frac{a}{2}\overline{|\nabla_x\Phi_{\overline{\rho}}|^2}\diff x, \qquad \cE_{rel}^{mv}(\t)=\int_{\O}\frac{1}{2}\overline{\rho|\bu-\bU|^2} \dx + \int_\Omega \frac{a}{2}\overline{|\nabla_x\Phi_{\overline{\rho}} - \nabla_x\Phi_r|^2}\diff x.
\end{align}
To avoid confusion we further denote the strong energy by
\begin{align}
    \cE^s(\tau)=\int_{\O}\frac{1}{2}r|U|^2 \dx + \int_\Omega\frac{a}{2}|\nabla_x\Phi_r|^2\diff x.
\end{align}

\begin{theorem}\label{t:main}
    Let $(\nu, \lambda, \nu^\infty)$ be a measure-valued solution in the sense of Definition \ref{def:measure_solutions} and $(r,\bU)\in C^{1}([0, T]\times \Omega)\times (C^{1}([0, T]\times \Omega)\cap L^\infty(0, T; W^{1,\infty}(\Omega)))$ be a strong solution with $r>0$ (and such that $r\to 0$, as $|x|\to +\infty$ in case $\Omega = \R^d$) to the pressureless Euler equations \eqref{pESdiv} as well as $\Psi$ satisfy the Assumption \ref{ass:on_kernel_psi}. Let $\e\geq0$ and suppose further that
 \begin{align*}
     \cE^{mv}_{rel}(0) + \|\rho_0 - r_0\|_{TV} + c(\O)\||x|^2(\rho_0 - r_0)\|_{TV}\leq \e,
 \end{align*}
 Then the relative entropy satisfies
 \begin{align*}
\cE_{rel}^{mv}(\t)\leq \e C_1e^{C_2\tau}
 \end{align*}
 for some constants $C_1,C_2>0$. Moreover, if
 \begin{align}
     \|\overline{\rho}_0-r_0\|_{TV}+\|\overline{\rho \bu}_0-r\bU_0\|_{TV}+\left|\cE^s(0)-\cE^{mv}(0)\right|=0,
 \end{align} then the following equalities hold almost everywhere in $(0,T) \times \O$:
 \begin{align*}
    \overline{\rho|\bu|^2}&=\overline{\rho}|\bU|^2,\\
    \overline{\rho\bu}&=\overline{\rho}\bU,\\
    \overline{\rho\bu\otimes\bu}&=\overline{\rho}\bU\otimes\bU,\\
    \overline{\rho}&=r.
 \end{align*}
\end{theorem}

\begin{remark}
    We refer to this theorem as weak-strong stability and uniqueness for the measure-valued solution, so that even if initial data does lead to finite time breakdown of the strong solution, the measure-valued solution remains close to strong solutions emanating from nearby initial data which may have a greater existence time.
\end{remark}

\section{Figalli-Kang type estimate for measure-valued solutions to continuity equations}\label{S:3}

In the present section, we shall prove a theorem reminiscent of the one given by Figalli and Kang related to the measure solutions to the continuity equations, cf. \cite{FK19}, where they have shown that a distance between a weak density and a strong density in the sense of Wasserstein can be estimated by the relative kinetic energy. In spirit, their results says the following: take a narrowly continuous measure solutions $\mu_t$ to the continuity equation
$$
\partial_t\mu_t + \diver_x(\mu_t\bu) = 0,
$$
for a velocity field $\bu(t, \cdot)\in L^p(\Omega, \diff\mu_t)$ for $p > 1$ and initial data $\mu_0$, and a strong solution $\rho$ of the continuity equation
$$
\partial_t\rho + \diver_x(\rho\bU) = 0,
$$
for a Lipschitz velocity field $\bU$, with initial data $\rho_0$. Then, one can control a $2$-Wasserstein distance between them by the inequality
$$
W_2^2(\mu_t, \rho_t) \leq C\int_{\Omega}|\bu_t - \bU_t|^2\diff\mu_t(x) + C\|\mu_0 - \rho_0\|_{TV} + Cc(\O)\||x|^2(\mu_0 - \rho_0)\|_{TV}.
$$
 We would like to have a similar inequality but in the sense of measure-valued solutions generated by generalized Young measures. We give a precise statement below.

\begin{lemma}\label{ineqFK}
    Let 
    $$
    \overline{\rho} = \langle \nu_{t,x},s^2\rangle+\langle \nu_{t,x}^\infty,s^2_\infty\rangle\l,
    $$
be the measure-valued density, which solves \eqref{CE-mv}, for the triple $(\nu, \lambda, \nu^\infty)$ of Definition \ref{def:measure_solutions}. Further suppose that $r$ solves
    $$
    \partial_t r + \diver_x(r\bU) = 0,
    $$
    with initial condition $r_0\in\mathcal{P}(\Omega)$ and $\bU\in L^\infty(0, T; W^{1, \infty}(\Omega))$. Then, there exists a constant $C > 0$ such that
    \begin{align*}
        W_2^2&(\overline{\rho}, r) \leq C\int_{\O}\overline{\rho|\bu-\bU|^2} \diff x + C\|\overline{\rho}_0 - r_0\|_{TV} + Cc(\O)\||x|^2(\overline{\rho}_0 - r_0)\|_{TV}\\
        &=C\left(\int_{\Omega}\langle\nu_{t, x}, |\beta - s\bU|^2\rangle\diff x + \int_{\Omega}\langle \nu_{t, x}^\infty, |\beta_\infty - s_\infty \bU|^2\rangle\diff\lambda_t(x)\right) + C\|\overline{\rho}_0 - r_0\|_{TV} + Cc(\O)\||x|^2(\overline{\rho}_0 - r_0)\|_{TV}.
    \end{align*}
\end{lemma}

\begin{proof}
    For the clarity of presentation let us define
    $$
    \diff\rho_1 := \langle \nu_{t,x},s^2\rangle\diff x,\qquad \diff\rho_2:=\langle \nu_{t,x}^\infty,s^2_\infty\rangle\diff\lambda,
    $$
    as well as
    $$
    \bu_1 := \frac{\langle \nu_{t,x},s\,\beta\rangle}{\langle \nu_{t,x},s^2\rangle},\qquad \bu_2 := \frac{\langle \nu_{t,x}^\infty,s_\infty\beta_\infty\rangle}{\langle\nu_{t,x}^\infty, s_\infty^2\rangle}.
    $$
    Then, the equation \eqref{CE-mv} can be written equivalently as
    \begin{align}\label{eq:rho_1_rho_2_form}
    \partial_t(\rho_1 + \rho_2) + \diver_x(\rho_1\bu_1 + \rho_2\bu_2) = 0,
    \end{align}
    in the sense of operators on $C^1([0, T]\times\overline{\Omega})$. Now notice that $\rho_1$ and $\rho_2$ are absolutely continuous with respect to the measure $\rho_1 + \rho_2$, hence due to the Radon--Nikodym theorem we can conclude that there exist Borel measurable $f_1$ and $f_2$ in $L^1(\Omega, \diff(\rho_1 +\rho_2))$ such that
    $$
    \diff\rho_1 = f_1\diff(\rho_1 + \rho_2), \qquad \diff\rho_2 = f_2\diff(\rho_1 + \rho_2).
    $$
    Moreover, it is easy to verify that $f_1$ and $f_2$ satisfy
    $$
    f_1, f_2 \geq 0,\qquad f_1 + f_2 = 1,\,(\rho_1 + \rho_2)-\text{ almost everywhere}.
    $$
    Thus, \eqref{eq:rho_1_rho_2_form} can be reformulated as
    $$
    \partial_t(\rho_1 + \rho_2) + \diver_x((\rho_1 + \rho_2)(f_1\bu_1 + f_2\bu_2)) = 0,
    $$
    in the sense of operators on $C^1([0, T]\times\overline{\Omega})$. Now, we will verify that our formulation is reasonable in some sense (cf. \cite[Proposition 3.3]{FK19}), which means that the velocity is in $L^p$ space with respect to $\rho_1 + \rho_2$ measure (in our case in $L^2$). Thanks to the Jensen inequality, and then the H\"{o}lder inequality, we may estimate
    \begin{align*}
        &\int_{\Omega}|f_1\bu_1 + f_2\bu_2|^2\diff(\rho_1 + \rho_2)\\
        &\leq \int_{\Omega}f_1|\bu_1|^2 + f_2|\bu_2|^2\diff(\rho_1 + \rho_2) = \int_{\Omega}|\bu_1|^2\diff\rho_1 + \int_{\Omega}|\bu_2|^2\diff\rho_2\\
        &\quad= \int_{\Omega}\left|\frac{\langle \nu_{t,x},s\,\beta\rangle}{\langle \nu_{t,x},s^2\rangle}\right|^2\langle\nu_{t,x}, s^2\rangle\diff x + \int_{\Omega}\left|\frac{\langle \nu_{t,x}^\infty,s_\infty\beta_\infty\rangle}{\langle\nu_{t,x}^\infty, s_\infty^2\rangle}\right|^2\langle\nu_{t,x}s_\infty^2\rangle\diff\lambda_t(x)\\
        &\leq \int_{\Omega}\langle\nu_{t,x},|\beta|^2\rangle\diff x + \int_{\Omega}\langle\nu_{t,x}, |\beta_\infty|^2\rangle\diff\lambda_t(x),
    \end{align*}
    which is finite due to the energy inequality \eqref{EN-mv}. Hence, we may apply the theory of Figalli and Kang \cite{FK19} to deduce that there exists a constant $C> 0$, such that
    \begin{align*}
        W_2^2(\rho_1 + \rho_2, r) \leq C\int_{\Omega}|f_1\bu_1 + f_2\bu_2 - \bU|^2\diff(\rho_1 + \rho_2) + C\|\overline{\rho}_0 - r_0\|_{TV} + Cc(\O)\||x|^2(\overline{\rho}_0 - r_0)\|_{TV}.
    \end{align*}
    We will now verify that
    $$
    \int_{\Omega}|f_1\bu_1 + f_2\bu_2 - \bU|^2\diff(\rho_1 + \rho_2) \leq \int_{\Omega}\langle\nu_{t, x}, |\beta - s\bU|^2\rangle\diff x + \int_{\Omega}\langle \nu_{t, x}^\infty, |\beta_\infty - s_\infty \bU|^2\rangle\diff\lambda_t(x),
    $$
    which will end the proof. Indeed, using the Jensen inequality, and then the H\"{o}lder inequality we may write
    \begin{align*}
        &\int_{\Omega}|f_1\bu_1 + f_2\bu_2 - \bU|^2\diff(\rho_1 + \rho_2) = \int_{\Omega}|f_1(\bu_1 - \bU) + f_2(\bu_2 - \bU)|^2\diff(\rho_1 + \rho_2)\\
        &\leq \int_{\Omega}f_1|\bu_1 - \bU|^2\diff(\rho_1 + \rho_2) + \int_{\Omega}f_2|\bu_2 - \bU|^2\diff(\rho_1 + \rho_2)\\
        &\quad = \int_{\Omega}|\bu_1 - \bU|^2\diff\rho_1 + \int_{\Omega}|\bu_2 - \bU|^2\diff\rho_2\\
        &\quad = \int_{\Omega}\frac{1}{\langle\nu_{t,x}, s^2\rangle}|\langle\nu_{t,x}, s(\beta - s\bU)\rangle|^2\diff x + \int_{\Omega}\frac{1}{\langle\nu_{t,x}, s_\infty^2\rangle}|\langle\nu^\infty_{t,x}, s_\infty(\beta_\infty - s_\infty\bU)\rangle|^2\diff \lambda_t(x)\\
        &\leq \int_{\Omega}\langle\nu_{t, x}, |\beta - s\bU|^2\rangle\diff x + \int_{\Omega}\langle \nu_{t, x}^\infty, |\beta_\infty - s_\infty \bU|^2\rangle\diff\lambda_t(x).
    \end{align*}
\end{proof}

\section{Relative Entropy Inequality}\label{S:4}
Within this section we establish a relative entropy inequality which provides stability of measure-valued solutions when the initial data is close to that of a strong solution. This is also a key step in proving the weak-strong stability and uniqueness result, Theorem \ref{t:main}.

It should be noted that without pressure the energy is not strictly convex in $\rho$ and hence the standard relative entropy inequality will fail to control what happens to the velocity in regions of vacuum, $\rho=0$. However, we can still obtain vanishing of the kinetic part of the relative entropy, and the Figalli-Kang type estimate provided in Lemma \ref{ineqFK} will provide us the relevant information on the densities.

We recall the definitions of the measure-valued energy and relative entropy
\begin{align*}
    \cE^{mv}(\tau)=\int_\O\frac{1}{2}\overline{\rho|\bu|^2} \dx + \int_\Omega \frac{a}{2}\overline{|\nabla_x\Phi_{\overline{\rho}}|^2}\diff x, \qquad \cE_{rel}^{mv}(\t)=\int_{\O}\frac{1}{2}\overline{\rho|\bu-\bU|^2} \dx + \int_\Omega \frac{a}{2}\overline{|\nabla_x\Phi_{\overline{\rho}} - \nabla_x\Phi_r|^2}\diff x.
\end{align*}
Note that the relative kinetic energy is the first term in the Figalli-Kang inequality, so that
\begin{align}
    W_2^2(\overline{\rho},r)\leq C\cE_{rel}^{mv}(\t)+C\|\overline{\rho}_0-r_0\|_{TV} + Cc(\O)\||x|^2(\overline{\rho}_0 - r_0)\|_{TV}.
\end{align}

\begin{lemma}\label{l:RelEnt}
    Let $(\nu, \lambda, \nu^\infty)$ be a measure-valued solution of \eqref{pESdiv} in the sense of Definition \ref{def:measure_solutions} and $(r,\bU)\in C^{1}([0, T]\times \Omega)\times (C^{1}([0, T]\times \Omega)\cap L^\infty(0, T; W^{1,\infty}(\Omega)))$ be a strong solution with $r>0$ (and such that $r\to 0$, as $|x|\to +\infty$ in case $\Omega = \R^d$) to the pressureless Euler equations \eqref{pESdiv} as well as $\Psi$ satisfy the Assumption \ref{ass:on_kernel_psi}. Let $\e\geq0$ and suppose further that
 \begin{align}
     \cE^{mv}_{rel}(0) + \|\rho_0 - r_0\|_{TV} + c(\O)\||x|^2(\rho_0 - r_0)\|_{TV}\leq \e
 \end{align}
 Then the relative entropy satisfies
 \begin{align}
\cE_{rel}^{mv}(\t)\leq \e C_1e^{C_2\tau}
 \end{align}
 for some constants $C_1,C_2>0$.
 \end{lemma}
 \begin{proof}
     We begin by expanding the relative entropy. 
     \begin{equation}\label{eq:relative_tau_expanded}
     \begin{split}
         \cE_{rel}^{mv}(\t)&=\int_{\O}\frac{1}{2}\overline{\rho|\bu|^2}\dx+\int_{\O}\frac{1}{2}\overline{\rho}|\bU|^2\dx-\int_{\O}\overline{\rho\bu}\cdot\bU\dx,\\
         &\quad +\int_{\Omega}\frac{a}{2}\overline{|\nabla_x\Phi_{\overline{\rho}}|^2}\diff x + \int_\Omega\frac{a}{2}|\nabla_x\Phi_r|^2\diff x - a\int_{\Omega}\overline{\nabla_x\Phi_{\overline{\rho}}}\cdot\nabla_x\Phi_r\diff x \\
         &=\cE^{mv}(\tau)+\int_{\O}\frac{1}{2}\overline{\rho}|\bU|^2\dx-\int_{\O}\overline{\rho\bu}\cdot\bU\dx\\
         &\quad + \int_\Omega\frac{a}{2}|\nabla_x\Phi_r|^2\diff x - a\int_{\Omega}\overline{\nabla_x\Phi_{\overline{\rho}}}\cdot\nabla_x\Phi_r\diff x.
        \end{split}
     \end{equation}
     If we test the continuity equation with $\frac{1}{2}|\bU|^2$ and $a\Phi_r$, and the momentum equation with $\bU$ we get
     \begin{equation}\label{eq:equations_testing_for_relative}
     \begin{split}
         \int_{\O}\frac{1}{2}\overline{\rho}|\bU|^2(\t,x)\dx&=\int_{\O}\frac{1}{2}\overline{\rho}|\bU|^2(0,x)+\int_0^\t\int_\O \overline{\rho}\bU\cdot \p_t\bU+\overline{\rho \bu}\cdot \n_x\bU\bU \dx\dt,\\
         a\int_{\O}\overline{\rho}\,\Phi_r(\t,x)\dx&=a\int_{\O}\overline{\rho}\,\Phi_r(0,x)+a\int_0^\t\int_\O \overline{\rho}\, \p_t\Phi_r+\overline{\rho \bu}\cdot \n_x\Phi_r \dx\dt,\\
         \int_{\O}\overline{\rho\bu}\cdot\bU(\t,x)\dx&=\int_{\O} (\overline{\rho\bu}\cdot \bU)(0,x) \dx+\int_0^\t\int_\O \overline{\rho\bu}\cdot \p_t\bU +(\overline{\rho\bu \otimes \bu}): \n_x\bU \dx \dt \\
         &+ \int_{0}^{\tau}\int_{\Omega}^{}\frac{a}{2}\overline{|\nabla_x\Phi_{\overline{\rho}}|^2}\operatorname{div}_x\bU\diff x\diff t + \int_0^\tau\int_{\Omega}a\overline{\nabla_x\Phi_{\overline{\rho}}\otimes\nabla_x\Phi_{\overline{\rho}}}:\nabla_x\bU\diff x\diff t\\
         &- \int_0^\tau\int_{\Omega}a M_{\overline{\rho}}\overline{\nabla_x\Phi_{\overline{\rho}}}\cdot\bU\,\diff x\diff t\\
         &-\frac{b}{2}\int_0^\tau\int_{\Omega\times\Omega}\overline{\rho(x)\rho(y)(\bu(x) - \bu(y))}\Psi(x-y)(\bU(x) - \bU(y))\diff y\diff x\diff t
        \end{split}
     \end{equation}
     Here, note that using the equation \eqref{eq:poissonmvs}
     \begin{align*}
         a\int_{\O}\overline{\rho}\,\Phi_r(\t,x)\dx = a\int_{\Omega}\overline{\nabla_x\Phi_{\overline{\rho}}}\cdot\nabla_x\Phi_r.
     \end{align*}
     Hence, as 
     $$
     \mathcal{E}^{mv}(\tau)\leq \cE^{mv}(0) - \frac{b}{2}\int_0^\tau\int_{\Omega\times \Omega}\overline{\rho(x)\rho(y)(\bu(x) - \bu(y))\Psi(x-y)(\bu(x) - \bu(y))}\diff y\diff x\diff t,
     $$
     we can combine \eqref{eq:relative_tau_expanded} with \eqref{eq:equations_testing_for_relative} to obtain
     \begin{equation}\label{ineq:relative_entropy_1}
     \begin{split}
         \cE_{rel}^{mv}(\t) - \cE^{mv}_{rel}(0)& +\frac{b}{2}\int_0^\tau\int_{\Omega\times \Omega}\overline{\rho(x)\rho(y)(\bu(x) - \bu(y))\Psi(x-y)(\bu(x) - \bu(y))}\diff y\diff x\diff t\\
         &\leq \int_\Omega\frac{a}{2}|\nabla_x\Phi_r|^2(\tau, x)\diff x -\int_\Omega\frac{a}{2}|\nabla_x\Phi_r|^2(0,x)\diff x\\
         &\qquad+ \int_0^\t\int_\O \overline{\rho}\bU\cdot \p_t\bU +\overline{\rho \bu}\cdot \n_x\bU\bU \dx\dt - \int_0^\t\int_\O \overline{\rho\bu}\cdot \p_t\bU +(\overline{\rho\bu \otimes \bu}): \n_x\bU \dx \dt\\
         &\qquad - a\int_0^\t\int_\O \overline{\rho}\, \p_t\Phi_r+\overline{\rho \bu}\cdot \n_x\Phi_r \dx\dt - \int_{0}^{\tau}\int_{\Omega}^{}\frac{a}{2}\overline{|\nabla_x\Phi_{\overline{\rho}}|^2}\operatorname{div}_x\bU\diff x\diff t \\
         &\qquad- \int_0^\tau\int_{\Omega}a\overline{\nabla_x\Phi_{\overline{\rho}}\otimes\nabla_x\Phi_{\overline{\rho}}}:\nabla_x\bU\diff x\diff t + \int_0^\tau\int_{\Omega}a M_{\overline{\rho}}\overline{\nabla_x\Phi_{\overline{\rho}}}\cdot\bU\,\diff x\diff t\\
         &\qquad +\frac{b}{2}\int_0^\tau\int_{\Omega\times\Omega}\overline{\rho(x)\rho(y)(\bu(x) - \bu(y))}\Psi(x-y)(\bU(x) - \bU(y))\diff y\diff x\diff t.
     \end{split}
     \end{equation}
     Now focusing on the terms without the potential $\Phi_r$ and non-local kernel $\Psi$, using the fact that for the strong solution with $r>0$ provides $\p_t\bU+ \n\bU \bU = -a\nabla_x\Phi_r + b\,\Psi\star(r\bU) - b\,\bU\Psi\star r$, we get
     \begin{align*}
        &\int_0^\t\int_\O \overline{\rho}\bU\cdot \p_t\bU +\overline{\rho \bu}\cdot \n_x\bU\bU \dx\dt - \int_0^\t\int_\O \overline{\rho\bu}\cdot \p_t\bU +(\overline{\rho\bu \otimes \bu}): \n_x\bU \dx \dt\\
        &\quad = \int_0^\tau\int_\Omega\overline{\rho\, (\bU - \bu)}\cdot\partial_t\bU + \overline{\rho(\bU - \bu)}\cdot\nabla_x\bU\bU + \overline{\rho(\bu - \bU)\otimes (\bU - \bu)}:\nabla_x\bU\diff x\diff t\\
        &\quad = -\int_0^\tau\int_\Omega a\overline{\rho(\bU - \bu)}\cdot\nabla_x\Phi_r - \overline{\rho(\bu - \bU)\otimes (\bU - \bu)}:\nabla_x\bU\diff x\diff t\\
        &\quad\qquad +b\int_0^\tau\int_{\Omega\times\Omega}\overline{\rho(\bU - \bu)}(x)\Psi(x-y)(r(y)\bU(y) - r(y)\bU(x))\diff y\diff x\diff t.
     \end{align*}
     Thus, \eqref{ineq:relative_entropy_1} reduces to
     \begin{equation}\label{ineq:relative_entropy_2}
     \begin{split}
         \cE_{rel}^{mv}(\t) - \cE^{mv}_{rel}(0)& +\frac{b}{2}\int_0^\tau\int_{\Omega\times \Omega}\overline{\rho(x)\rho(y)(\bu(x) - \bu(y))\Psi(x-y)(\bu(x) - \bu(y))}\diff y\diff x\diff t\\
         &\leq \int_\Omega\frac{a}{2}|\nabla_x\Phi_r|^2(\tau, x)\diff x -\int_\Omega\frac{a}{2}|\nabla_x\Phi_r|^2(0,x)\diff x\\
         &\qquad- \int_0^\tau\int_\Omega a\overline{\rho(\bU - \bu)}\cdot\nabla_x\Phi_r - \overline{\rho(\bu - \bU)\otimes (\bU - \bu)}:\nabla_x\bU\diff x\diff t\\
         &\qquad - a\int_0^\t\int_\O \overline{\rho}\, \p_t\Phi_r+\overline{\rho \bu}\cdot \n_x\Phi_r \dx\dt - \int_{0}^{\tau}\int_{\Omega}^{}\frac{a}{2}\overline{|\nabla_x\Phi_{\overline{\rho}}|^2}\operatorname{div}_x\bU\diff x\diff t \\
         &\qquad- \int_0^\tau\int_{\Omega}a\overline{\nabla_x\Phi_{\overline{\rho}}\otimes\nabla_x\Phi_{\overline{\rho}}}:\nabla_x\bU\diff x\diff t + \int_0^\tau\int_{\Omega}a M_{\overline{\rho}}\overline{\nabla_x\Phi_{\overline{\rho}}}\cdot\bU\,\diff x\diff t\\
         &\qquad +\frac{b}{2}\int_0^\tau\int_{\Omega\times\Omega}\overline{\rho(x)\rho(y)(\bu(x) - \bu(y))}\Psi(x-y)(\bU(x) - \bU(y))\diff y\diff x\diff t\\
         &\qquad +b\int_0^\tau\int_{\Omega\times\Omega}\overline{\rho(\bU - \bu)}(x)\Psi(x-y)(r(y)\bU(y) - r(y)\bU(x))\diff y\diff x\diff t.
     \end{split}
     \end{equation}
     We turn our attention to the terms with the potential $\Phi_r$. Using the Poisson's equation for $\Phi_r$ we get
     \begin{equation}\label{nonlocal_reduction_1}
         \begin{split}
             &\int_\Omega\frac{a}{2}|\nabla_x\Phi_r|^2(\tau, x)\diff x -\int_\Omega\frac{a}{2}|\nabla_x\Phi_r|^2(0,x)\diff x -a\int_0^\tau\int_\Omega \overline{\rho(\bU - \bu)}\cdot\nabla_x\Phi_r\diff x\diff t\\
             &\qquad- a\int_0^\t\int_\O \overline{\rho}\, \p_t\Phi_r+\overline{\rho \bu}\cdot \n_x\Phi_r \dx\dt - \int_{0}^{\tau}\int_{\Omega}^{}\frac{a}{2}\overline{|\nabla_x\Phi_{\overline{\rho}}|^2}\operatorname{div}_x\bU\diff x\diff t \\
         &\qquad- \int_0^\tau\int_{\Omega}a\overline{\nabla_x\Phi_{\overline{\rho}}\otimes\nabla_x\Phi_{\overline{\rho}}}:\nabla_x\bU\diff x\diff t + \int_0^\tau\int_{\Omega}a M_{\overline{\rho}}\overline{\nabla_x\Phi_{\overline{\rho}}}\cdot\bU\,\diff x\diff t.\\
             & = \int_\Omega\frac{a}{2}|\nabla_x\Phi_r|^2(\tau, x)\diff x -\int_\Omega\frac{a}{2}|\nabla_x\Phi_r|^2(0,x)\diff x - a\int_0^\t\int_\O \overline{\rho}\, \p_t\Phi_r+\overline{\rho} \bU\cdot \n_x\Phi_r \dx\dt\\
             &\qquad - \int_{0}^{\tau}\int_{\Omega}^{}\frac{a}{2}\overline{|\nabla_x\Phi_{\overline{\rho}}|^2}\operatorname{div}_x\bU\diff x\diff t - \int_0^\tau\int_{\Omega}a\overline{\nabla_x\Phi_{\overline{\rho}}\otimes\nabla_x\Phi_{\overline{\rho}}}:\nabla_x\bU\diff x\diff t + \int_0^\tau\int_{\Omega}a M_{\overline{\rho}}\overline{\nabla_x\Phi_{\overline{\rho}}}\cdot\bU\,\diff x\diff t\\
             & = a\int_0^\tau\int_\Omega \nabla_x\Phi_r\cdot\partial_t\nabla_x\Phi_r\diff x\diff t - a\int_0^\t\int_\O \overline{\rho}\, \p_t\Phi_r+\overline{\rho} \bU\cdot \n_x\Phi_r \dx\dt\\
             &\qquad - \int_{0}^{\tau}\int_{\Omega}^{}\frac{a}{2}\overline{|\nabla_x\Phi_{\overline{\rho}}|^2}\operatorname{div}_x\bU\diff x\diff t - \int_0^\tau\int_{\Omega}a\overline{\nabla_x\Phi_{\overline{\rho}}\otimes\nabla_x\Phi_{\overline{\rho}}}:\nabla_x\bU\diff x\diff t + \int_0^\tau\int_{\Omega}a M_{\overline{\rho}}\overline{\nabla_x\Phi_{\overline{\rho}}}\cdot\bU\,\diff x\diff t\\
             & = a\int_0^\tau\int_\Omega r\partial_t \Phi_r\diff x\diff t- a\int_0^\t\int_\O \overline{\rho}\, \p_t\Phi_r+\overline{\rho} \bU\cdot \n_x\Phi_r \dx\dt - a\int_0^\tau\int_\Omega M_r\partial_t \Phi_r\diff x\diff t\\
             &\qquad - \int_{0}^{\tau}\int_{\Omega}^{}\frac{a}{2}\overline{|\nabla_x\Phi_{\overline{\rho}}|^2}\operatorname{div}_x\bU\diff x\diff t - \int_0^\tau\int_{\Omega}a\overline{\nabla_x\Phi_{\overline{\rho}}\otimes\nabla_x\Phi_{\overline{\rho}}}:\nabla_x\bU\diff x\diff t + \int_0^\tau\int_{\Omega}a M_{\overline{\rho}}\overline{\nabla_x\Phi_{\overline{\rho}}}\cdot\bU\,\diff x\diff t\\
             & = a\int_0^\tau\int_\Omega (\nabla_x\Phi_r - \nabla_x\Phi_{\overline{\rho}})\cdot\nabla_x\partial_t\Phi_r\diff x\diff t - a\int_0^\t\int_\O \overline{\rho}\bU\cdot \n_x\Phi_r \dx\dt \\
             &\qquad - \int_{0}^{\tau}\int_{\Omega}^{}\frac{a}{2}\overline{|\nabla_x\Phi_{\overline{\rho}}|^2}\operatorname{div}_x\bU\diff x\diff t - \int_0^\tau\int_{\Omega}a\overline{\nabla_x\Phi_{\overline{\rho}}\otimes\nabla_x\Phi_{\overline{\rho}}}:\nabla_x\bU\diff x\diff t + \int_0^\tau\int_{\Omega}a M_{\overline{\rho}}\overline{\nabla_x\Phi_{\overline{\rho}}}\cdot\bU\,\diff x\diff t
         \end{split}
     \end{equation}
Before moving forward, let us notice that since $(r, \Phi_r)$ is a strong solution to $-\Delta\Phi_r  = r - M_r$ and $(r, \bU)$ is a strong solution to the continuity equation, then using the distributional formulation \eqref{eq:poissonmvs} for $(\overline{\rho}, \overline{\nabla_x\Phi_{\overline{\rho}}})$, we get
    \begin{equation*}
    \begin{split}
            &\int_0^\tau\int_{\Omega}(r - M_r)\overline{\nabla_x\Phi_{\overline{\rho}}}\cdot \bU \diff x\diff t + \int_0^\tau\int_{\Omega}(\overline{\rho} - M_{\overline{\rho}})\nabla_x\Phi_r\cdot \bU\diff x\diff t\\
            &=-\int_0^\tau\int_{\Omega}\Delta_x\Phi_r\,\overline{\nabla_x\Phi_{\overline{\rho}}}\cdot \diff x\diff t + \int_0^\tau\int_{\Omega}\nabla_x\Phi_{\overline{\rho}}\cdot\nabla_x(\nabla_x\Phi_r\cdot \bU)\diff x\diff t\\
            &= -\int_0^\tau\int_{\Omega}\Delta_x\Phi_r\,\overline{\nabla_x\Phi_{\overline{\rho}}}\cdot \bU \diff x\diff t + \int_0^\tau\int_{\Omega}\nabla_x\Phi_r\otimes\nabla_x\Phi_{\overline{\rho}}:\nabla_x \bU\diff x\diff t\\
            &\quad+ \int_0^\tau\int_{\Omega}\nabla_x\Phi_{\overline{\rho}}\cdot D^2_x\Phi_r\cdot \bU\diff x\diff t \\
            &  =\int_0^\tau\int_\Omega\sum_{i, j = 1}^d\frac{\p}{\p x_j}\left(\frac{\p\Phi_r}{\p x_i}U_j - \frac{\p\Phi_r}{\p x_j}U_i\right)\frac{\p\Phi_{\overline{\rho}}}{\p x_i}\diff x\diff t -\int_0^\tau\int_{\Omega}\overline{\nabla_x\Phi_{\overline{\rho}}}\cdot\nabla_x\Phi_r\,\diver_x \bU\diff x\diff t \\
            &\quad+ \int_0^\tau\int_{\Omega}\overline{\nabla_x\Phi_{\overline{\rho}}}\otimes\nabla_x\Phi_r : \nabla_x \bU\diff x\diff t + \int_0^\tau\int_{\Omega}\nabla_x\Phi_{r}\otimes\overline{\nabla_x\Phi_{\overline{\rho}}} : \nabla_x \bU\diff x\diff t,
    \end{split}
    \end{equation*}
which, since by integration by parts
\begin{equation*}
    \begin{split}
        &\int_0^\tau\int_\Omega\sum_{i, j = 1}^d\frac{\p}{\p x_j}\left(\frac{\p\Phi_r}{\p x_i}U_j - \frac{\p\Phi_r}{\p x_j}U_i\right)\frac{\p\Phi_{\overline{\rho}}}{\p x_i}\diff x\diff t = \int_0^\tau\int_\Omega\sum_{i, j = 1}^d\frac{\p}{\p x_i}\left(\frac{\p\Phi_r}{\p x_i}U_j - \frac{\p\Phi_r}{\p x_j}U_i\right)\frac{\p\Phi_{\overline{\rho}}}{\p x_j}\diff x\diff t\\
        & = -\int_0^\tau\int_\Omega\sum_{i, j = 1}^d\frac{\p}{\p x_j}\left(\frac{\p\Phi_r}{\p x_i}U_j - \frac{\p\Phi_r}{\p x_j}U_i\right)\frac{\p\Phi_{\overline{\rho}}}{\p x_i}\diff x\diff t,
    \end{split}
\end{equation*}
implies
\begin{equation}\label{eq:reductionpoisson4}
    \begin{split}
        &\int_0^\tau\int_{\Omega}(r - M_r)\overline{\nabla_x\Phi_{\overline{\rho}}}\cdot \bU \diff x\diff t + \int_0^\tau\int_{\Omega}(\overline{\rho} - M_{\overline{\rho}})\nabla_x\Phi_r\cdot \bU\diff x\diff t\\
        &  = -\int_0^\tau\int_{\Omega}\overline{\nabla_x\Phi_{\overline{\rho}}}\cdot\nabla_x\Phi_r\,\diver_x \bU\diff x\diff t \\
            &\quad+ \int_0^\tau\int_{\Omega}\overline{\nabla_x\Phi_{\overline{\rho}}}\otimes\nabla_x\Phi_r : \nabla_x \bU\diff x\diff t + \int_0^\tau\int_{\Omega}\nabla_x\Phi_{r}\otimes\overline{\nabla_x\Phi_{\overline{\rho}}} : \nabla_x \bU\diff x\diff t.
    \end{split}
\end{equation}
Using \eqref{eq:reductionpoisson4} we may now calculate
\begin{equation*}
    \begin{split}
             &- \int_{0}^{\tau}\int_{\Omega}^{}\frac{a}{2}\overline{|\nabla_x\Phi_{\overline{\rho}}|^2}\operatorname{div}_x\bU\diff x\diff t + a\overline{\nabla_x\Phi_{\overline{\rho}}\otimes\nabla_x\Phi_{\overline{\rho}}}:\nabla_x\bU\diff x\diff t- a M_{\overline{\rho}}\overline{\nabla_x\Phi_{\overline{\rho}}}\cdot\bU + a\overline{\rho}\bU\cdot \n_x\Phi_r\diff x\diff t\\
             & = -\frac{a}{2}\int_0^\tau\int_\Omega\overline{|\nabla_x\Phi_{\overline{\rho}} - \nabla_x\Phi_r|^2}\diver_x\bU\diff x\diff t + a\int_0^\tau\int_{\Omega}\overline{(\nabla_x\Phi_{\overline{\rho}} - \nabla_x\Phi_r)\otimes (\nabla_x\Phi_{\overline{\rho}} - \nabla_x\Phi_r)}: \nabla_x \bU\diff x\diff t \\
             &\qquad + a\int_0^\tau\int_\Omega \overline{(\nabla_x\Phi_{\overline{\rho}} - \nabla_x\Phi_r)}\cdot r\bU\diff x\diff t,
    \end{split}
\end{equation*}
which combined with \eqref{nonlocal_reduction_1} grants us
\begin{equation}\label{eq:reduction_poisson_5}
    \begin{split}
        &\int_\Omega\frac{a}{2}|\nabla_x\Phi_r|^2(\tau, x)\diff x -\int_\Omega\frac{a}{2}|\nabla_x\Phi_r|^2(0,x)\diff x -a\int_0^\tau\int_\Omega \overline{\rho(\bU - \bu)}\cdot\nabla_x\Phi_r\diff x\diff t\\
             &\qquad- a\int_0^\t\int_\O \overline{\rho}\, \p_t\Phi_r+\overline{\rho \bu}\cdot \n_x\Phi_r \dx\dt - \int_{0}^{\tau}\int_{\Omega}^{}\frac{a}{2}\overline{|\nabla_x\Phi_{\overline{\rho}}|^2}\operatorname{div}_x\bU\diff x\diff t \\
         &\qquad- \int_0^\tau\int_{\Omega}a\overline{\nabla_x\Phi_{\overline{\rho}}\otimes\nabla_x\Phi_{\overline{\rho}}}:\nabla_x\bU\diff x\diff t + \int_0^\tau\int_{\Omega}a M_{\overline{\rho}}\overline{\nabla_x\Phi_{\overline{\rho}}}\cdot\bU\,\diff x\diff t.\\
         & = -\frac{a}{2}\int_0^\tau\int_\Omega\overline{|\nabla_x\Phi_{\overline{\rho}} - \nabla_x\Phi_r|^2}\diver_x\bU\diff x\diff t + a\int_0^\tau\int_{\Omega}\overline{(\nabla_x\Phi_{\overline{\rho}} - \nabla_x\Phi_r)\otimes (\nabla_x\Phi_{\overline{\rho}} - \nabla_x\Phi_r)}: \nabla_x \bU\diff x\diff t \\
             &\qquad + a\int_0^\tau\int_\Omega (\nabla_x\Phi_r - \nabla_x\Phi_{\overline{\rho}})\cdot\nabla_x\partial_t\Phi_r\diff x\diff t + a\int_0^\tau\int_\Omega \overline{(\nabla_x\Phi_{\overline{\rho}} - \nabla_x\Phi_r)}\cdot r\bU\diff x\diff t
    \end{split}
\end{equation}
Applying \eqref{eq:poissonmvs}, the continuity equation for $r$ and the Poisson equation for $\Phi_r$ it is easy to verify, that
\begin{align*}
     a\int_0^\tau\int_\Omega (\nabla_x\Phi_r - \nabla_x\Phi_{\overline{\rho}})\cdot\nabla_x\partial_t\Phi_r\diff x\diff t + a\int_0^\tau\int_\Omega \overline{(\nabla_x\Phi_{\overline{\rho}} - \nabla_x\Phi_r)}\cdot r\bU\diff x\diff t = 0.
\end{align*}
Hence, putting \eqref{eq:reduction_poisson_5} into \eqref{ineq:relative_entropy_2} we obtain
\begin{equation}\label{ineq:relative_entropy_3}
     \begin{split}
         \cE_{rel}^{mv}(\t) - \cE^{mv}_{rel}(0)& +\frac{b}{2}\int_0^\tau\int_{\Omega\times \Omega}\overline{\rho(x)\rho(y)(\bu(x) - \bu(y))\Psi(x-y)(\bu(x) - \bu(y))}\diff y\diff x\diff t\\
         &\leq -\frac{a}{2}\int_0^\tau\int_\Omega\overline{|\nabla_x\Phi_{\overline{\rho}} - \nabla_x\Phi_r|^2}\diver_x\bU\diff x\diff t\\
         &\qquad + a\int_0^\tau\int_{\Omega}\overline{(\nabla_x\Phi_{\overline{\rho}} - \nabla_x\Phi_r)\otimes (\nabla_x\Phi_{\overline{\rho}} - \nabla_x\Phi_r)}: \nabla_x \bU\diff x\diff t\\
         &\qquad+\int_0^\tau\int_\Omega \overline{\rho(\bu - \bU)\otimes (\bU - \bu)}:\nabla_x\bU\diff x\diff t\\
         &\qquad +\frac{b}{2}\int_0^\tau\int_{\Omega\times\Omega}\overline{\rho(x)\rho(y)(\bu(x) - \bu(y))}\Psi(x-y)(\bU(x) - \bU(y))\diff y\diff x\diff t\\
         &\qquad + b\int_0^\tau\int_{\Omega\times\Omega}\overline{\rho(\bU - \bu)}(x)\Psi(x-y)(r(y)\bU(y) - r(y)\bU(x))\diff y\diff x\diff t.
     \end{split}
     \end{equation}
At last, we handle the non-local terms with potential $\Psi$. By the symmetry of the tensorized measure and the potential $\Psi$ we can write
\begin{align*}
\frac{b}{2}\int_0^\tau\int_{\Omega\times\Omega}\overline{\rho(x)\rho(y)(\bu(x) - \bu(y))}\Psi(x-y)(\bU(x) - \bU(y))\diff y\diff x\diff t \\
= b\int_0^\tau\int_{\Omega\times\Omega}\overline{\rho(x)\rho(y)(\bu(x) - \bu(y))}\Psi(x-y)\bU(x)\diff y\diff x\diff t.
\end{align*}
Now, we may split our two integrals in the following way
\begin{equation*}
    \begin{split}
    &- \int_0^\tau\int_{\Omega\times\Omega}r(y)\overline{\rho(\bU - \bu)}(x)\Psi(x-y)(\bU(x) - \bU(y))\diff y\diff x\diff t\\
    &+ \int_0^\tau\int_{\Omega\times\Omega}\overline{\rho(x)\rho(y)(\bu(x) - \bu(y))}\Psi(x-y)\bU(x)\diff y\diff x\diff t\\
    & \quad = -  \int_0^\tau\int_{\Omega\times\Omega}r(y)\overline{\rho(\bU - \bu)}(x)\Psi(x-y)(\bU(x) - \bU(y))\diff y\diff x\diff t\\
    &\qquad + \int_0^\tau\int_{\Omega\times\Omega}\overline{\rho(x)\rho(y)(\bU(x) - \bu(x))}\Psi(x-y)(\bU(x) - \bU(y))\diff y\diff x\diff t\\
    &\qquad - \int_0^\tau\int_{\Omega\times\Omega}\overline{\rho(x)\rho(y)(\bU(x) - \bu(x))}\Psi(x-y)(\bU(x) - \bU(y))\diff y\diff x\diff t\\
    &\qquad + \int_0^\tau\int_{\Omega\times\Omega}\overline{\rho(x)\rho(y)(\bU(x) - \bu(x))\Psi(x - y)(\bu(x) - \bu(y))}\diff y\diff x\diff t\\
    &\qquad + \int_0^\tau\int_{\Omega\times\Omega}\overline{\rho(x)\rho(y)\bu(x)\Psi(x - y)(\bu(x) - \bu( y))}\diff y \diff x\diff t\\
    &\quad =: I_1 + I_2 + I_3 +I_4 + I_5.
    \end{split}
\end{equation*}
Now, by Lemma \ref{lem:alignment_inequality_for_relative_entropy} we can estimate
\begin{equation}\label{ineq:entropy_I1_I2}
    \begin{split}
    & |I_1 + I_2|\\
    &=\Bigg|-\int_0^\tau\int_{\Omega\times\Omega}r(y)\overline{\rho(\bU - \bu)}(x)\Psi(x-y)(\bU(x) - \bU(y))\diff y\diff x\diff t\\
    &\qquad + \int_0^\tau\int_{\Omega\times\Omega}\overline{\rho(x)\rho(y)(\bU(x) - \bu(x))}\Psi(x-y)(\bU(x) - \bU(y))\diff y\diff x\diff t\Bigg|\\
    &= \Bigg|\int_0^\tau\int_{\Omega\times\Omega}\langle\nu^2_{t,x, y}, (s^2\bU(x) - s\beta)((s')^2 - r(y))\rangle\Psi(x-y)(\bU(x) - \bU(y))\diff y\diff x\diff t\\
    &\qquad + \int_0^\tau\int_{\Omega\times\Omega}\langle\nu^{2,\infty}_{t,x, y}, (s^2_\infty\bU(x) - s_\infty\beta_\infty)((s'_\infty)^2 - r(y))\rangle\Psi(x-y)(\bU(x) - \bU(y))\diff\l(y)\diff\l(x)\diff t\Bigg|\\
    &\leq C\|\Psi\|_{W^{1,\infty}}\|\bU\|_{L^\infty_tW^{1,\infty}_x}\left(\int_0^\tau\int_{\Omega}\overline{\rho|\bu - \bU|}^2\diff x\diff t + \int_0^\tau d_f^2(\overline{\rho}_t, r_t)\diff t\right)
    \end{split}
\end{equation}
Moving forward, using the fact that $\overline{\rho}_t$ is a probability measure, and the H\"{o}lder inequality 
\begin{equation}\label{ineq:entropy_I3_I4}
    \begin{split}
    &|I_3 + I_4|\\
    & = \left|\int_0^\tau\int_{\Omega\times\Omega}\overline{\rho(x)\rho(y)(\bU(x) - \bu(x))\Psi(x - y)(\bu(x) - \bU(x) - (\bu(y) - \bU(y)))}\diff y\diff x\diff t\right|\\
    & = \Bigg|\int_0^\tau\int_{\Omega\times\Omega}\langle\nu^2_{t, x, y}, (s')^2(s\bU(x) - \beta)\Psi(x-y)(\beta - s\bU(x))\rangle\diff x\diff y\diff t\\
    &\quad + \int_0^\tau\int_{\Omega\times\Omega}\langle\nu^2_{t, x, y}, s s'(s\bU(x) - \beta)\Psi(x-y)(\beta' - s'\bU(y))\rangle\diff x\diff y\diff t\\
    &\quad + \int_0^\tau\int_{\Omega\times\Omega}\langle\nu^{2,\infty}_{t, x, y}, (s_\infty')^2(s_\infty\bU(x) - \beta_\infty)\Psi(x-y)(\beta_\infty - s_\infty\bU(x))\rangle\diff\lambda_t(x)\diff\lambda_t(y)\diff t\\
    &\quad + \int_0^\tau\int_{\Omega\times\Omega}\langle\nu^{2,\infty}_{t, x, y}, s_\infty s_\infty'(s_\infty\bU(x) - \beta_\infty)\Psi(x-y)(\beta'_\infty - s'_\infty\bU(y))\rangle\diff \lambda_t(x)\diff \lambda_t(y)\diff t\Bigg|\\
    & \leq \|\Psi\|_\infty\Bigg|\int_0^\tau\int_{\Omega\times\Omega}\langle\nu_{t, x}, |s\bU(x) - \beta|^2\rangle\langle\nu_{t,y}, (s')^2\rangle\diff x\diff y\diff t\\
    &\quad + \int_0^\tau\int_{\Omega\times\Omega}\langle\nu^2_{t, x, y}, s s'|s\bU(x) - \beta||\beta' - s'\bU(y)|\rangle\diff x\diff y\diff t\\
    &\quad + \int_0^\tau\int_{\Omega\times\Omega}\langle\nu^{\infty}_{t, x},|s_\infty\bU(x) - \beta_\infty|^2\rangle\langle\nu^\infty_{t,y}, (s'_\infty)^2\rangle\diff\lambda_t(x)\diff\lambda_t(y)\diff t\\
    &\quad + \int_0^\tau\int_{\Omega\times\Omega}\langle\nu^{2,\infty}_{t, x, y}, s_\infty s_\infty'|s_\infty\bU(x) - \beta_\infty||\beta'_\infty - s'_\infty\bU(y)|\rangle\diff \lambda_t(x)\diff \lambda_t(y)\diff t\Bigg|\\
    & \leq \|\Psi\|_\infty\int_0^\tau\int_{\Omega}\overline{\rho|\bu - \bU|}^2\diff x\diff t + \|\Psi\|_\infty\int_0^\tau\left(\int_{\Omega}\overline{\rho|\bu - \bU|}\diff x\right)^2\diff t\\
    &\leq 2\|\Psi\|_\infty\int_0^\tau\int_{\Omega}\overline{\rho|\bu - \bU|}^2\diff x\diff t.
    \end{split}
\end{equation}
At last, the fact that $\Psi$ is an even function implies
$$
I_5 = \int_0^\tau\int_{\Omega\times\Omega}\overline{\rho(x)\rho(y)(\bu(x) - \bu(y))\Psi(x - y)(\bu(x) - \bu( y))}\diff y \diff x\diff t,
$$
which cancels out with the same term on the left-hand side of 
\eqref{ineq:relative_entropy_3}. Thus, putting \eqref{ineq:entropy_I1_I2} and \eqref{ineq:entropy_I3_I4} as well as utilizing standard inequalities grants us
\begin{equation}\label{ineq:relative_entropy_4}
    \begin{split}
        \cE_{rel}^{mv}(\t) - &\cE^{mv}_{rel}(0)\leq \|\bU\|_{L^\infty_t W^{1,\infty}_x}\frac{3a}{2}\int_0^\tau\int_\Omega\overline{|\nabla_x\Phi_{\overline{\rho}} - \nabla_x\Phi_r|^2}\diff x\diff t\\
         &\qquad+\left(\|\bU\|_{L^\infty_t W^{1,\infty}_x} + b\|\Psi\|_\infty + Cb\|\Psi\|_{W^{1,\infty}}\|\bU\|_{L^\infty_t W^{1,\infty}_x}\right)\int_0^\tau\int_\Omega \overline{\rho|\bu - \bU|}^2\diff x\diff t\\
         &\qquad + Cb\|\Psi\|_{W^{1,\infty}}\|\bU\|_{L^\infty_tW^{1,\infty}_x}\int_0^\tau d_f^2(\overline{\rho}_t, r_t)\diff t.
    \end{split}
\end{equation}
As at this point we lack the distance between densities in the relative entropy, we add the inequality from Lemma \ref{ineqFK} to \eqref{ineq:relative_entropy_4} and use the fact that the flat metric is bounded above by the $W_2$ distance to deduce
\begin{equation}\label{ineq:relative_entropy_5}
    \begin{split}
        \cE_{rel}^{mv}(\t) + W_2^2(\overline{\rho}_\t, r_\t)&\leq  \cE^{mv}_{rel}(0) +  \|\bU\|_{L^\infty_t W^{1,\infty}_x}\frac{3a}{2}\int_0^\tau\int_\Omega\overline{|\nabla_x\Phi_{\overline{\rho}} - \nabla_x\Phi_r|^2}\diff x\diff t\\
         &\qquad+\left(C +\|\bU\|_{L^\infty_t W^{1,\infty}_x} + b\|\Psi\|_\infty + Cb\|\Psi\|_{W^{1,\infty}}\|\bU\|_{L^\infty_t W^{1,\infty}_x}\right)\int_0^\tau\int_\Omega \overline{\rho|\bu - \bU|}^2\diff x\diff t\\
         &\qquad + Cb\|\Psi\|_{W^{1,\infty}}\|\bU\|_{L^\infty_tW^{1,\infty}_x}\int_0^\tau W_2^2(\overline{\rho}_t, r_t)\diff t\\
         &\qquad + C\|\rho_0 - r_0\|_{TV} + Cc(\O)\||x|^2(\rho_0 - r_0)\|_{TV}.
    \end{split}
\end{equation}
Application of the Gr\"{o}nwall inequality to \eqref{ineq:relative_entropy_5} grants us our thesis.
\end{proof}
This result establishes that measure-valued solutions remain close to corresponding strong solutions if initially so. In particular, Lemma \ref{l:RelEnt} yields the following corollary.
\begin{corollary}\label{c:relent}
    If $\e=0$ in the context of Lemma \ref{l:RelEnt} so that the measure-valued solution and the strong solution have the same initial data, then
    \begin{align*}
        \cE^{mv}_{rel}(\tau)=0, \quad \text{a.e} \ \t\in [0,T].
    \end{align*}
\end{corollary}

\section{Weak-Strong Uniqueness}\label{S:5}
For identical initial data Corollary \ref{c:relent} provides
\begin{align}
0=\int_{[0,T] \times \O} \frac{1}{2}\overline{\rho|\bu-\bU|^2} \dx \dt=
    \int_{[0,T]\times \O} \langle \nu_{t,x}, |\beta-s\bU|^2\rangle \dx \dt+\int_{[0,T]\times \O} \langle \nu_{t,x}^\infty, |\beta_\infty-s_\infty\bU|^2\rangle \dl_t(x) \dt
\end{align}

However, this does not immediately grant the weak-strong uniqueness property. We first state without proof several useful inequalities.
\begin{lemma}\label{usefulineq}
    Let $\bU \in L^\infty(\R^n;\R^d)$ be a bounded vector-valued function. Then for any vector $\bu \in \R^d$ and $\d>0$ we have the following inequalities:
    
         For any $i,j=1,...,d$,
         \begin{align}
         |u_iu_j-U_iU_j|&\leq c\d+C_\d(|u_i-U_i|^2+|u_j-U_j|^2),\\
         \left||\bu|^2-|\bU|^2\right|&\leq c\d+C_\d|\bu-\bU|^2,\\
         |\bu-\bU|&\leq c\d+C_\d|\bu-\bU|^2.
         \end{align}
   Where $c$ depends only on $\|\bU\|_{L^\infty}$ and $C_\d$ on $\d$ and $\|\bU\|_{L^\infty}$. 
\end{lemma}
We now seek to identify the measure-valued density $\overline{\rho}$ with the strong density $r$.

\begin{lemma}
    Let $(\nu, \lambda, \nu^\infty)$ be a measure-valued solution of \eqref{pESdiv} and $(r,\bU)\in C^{1}([0, T]\times \Omega)\times (C^{1}([0, T]\times \Omega)\cap L^\infty(0, T; W^{1,\infty}(\Omega)))$ be a strong solution with $r>0$ (and such that $r\to 0$, as $|x|\to +\infty$ in case $\Omega = \R^d$) to the pressureless Euler equations \eqref{pESdiv} as well as $\Psi$ satisfy the Assumption \ref{ass:on_kernel_psi}. Suppose further that
 \begin{align}
     \|\overline{\rho}_0-r_0\|_{TV}+\|\overline{\rho \bu}_0-r\bU_0\|_{TV}+\left|\cE^s(0)-\cE^{mv}(0)\right|=0,
 \end{align}
 then the following equalities hold almost everywhere in $(0,T) \times \O$:
 \begin{align}
    \overline{\rho|\bu|^2}&=\overline{\rho}|\bU|^2,\label{ident1}\\
    \overline{\rho\bu}&=\overline{\rho}\bU,\label{ident2}\\
    \overline{\rho\bu\otimes\bu}&=\overline{\rho}\bU\otimes\bU,\label{ident3}\\
    \overline{\rho}&=r.\label{ident4}
 \end{align}
\end{lemma}
\begin{proof}
    
We begin with \eqref{ident1}. Let $\d>0$.

Let $\phi \in C_c^{\infty}(\O)$ with $\phi\geq 0$ be a test function. Consider first by the triangle inequality,
\begin{align*}
    &\left|\int_{[0,T]\times\O}\phi\left(\overline{\rho|\bu|^2}-\overline{\rho}|\bU|^2\right)\dx\dt\right|\leq \int_{[0,T]\times\O}\phi\left|\overline{\rho|\bu|^2}-\overline{\rho}|\bU|^2\right|\dx\dt,\\
    &\leq\int_{[0,T]\times\O}\phi \left|\langle\nu_{t,x},(|\b|^2-s^2|\bU|^2) \rangle\right|\dx\dt+\int_{[0,T]\times\O}\phi \left|\langle\nu_{t,x}^\infty,(|\b_\infty|^2-s_{\infty}^2|\bU|^2) \rangle\right|\dl_t(x)\dt
\end{align*}
Then using Lemma \ref{usefulineq}
\begin{align*}
    &\leq \int_{[0,T]\times\O}\phi \left\langle\nu_{t,x},c\d s^2 +C_\d\left|\b-s\bU\right|^2 \right\rangle\dx \dt+\int_{[0,T]\times\O}\phi \left\langle\nu_{t,x}^\infty,c \d s^2_\infty+C_\d\left|\b_\infty-s_\infty\bU\right|^2 \right\rangle\dl_t(x) \dt,\\
    &=c\d \int_{[0,T]\times\O}\phi\overline{\rho} \dx \dt.
\end{align*}
Where we used the vanishing of the relative entropy. Now with conservation of mass, and the fact that $\delta>0$ was arbitrary, we can conclude
\begin{align*}
    \overline{\rho|\bu|^2}=\overline{\rho}|\bU|^2,
\end{align*}
almost everywhere in $(0,T)\times \Omega$.

To achieve 
\begin{align*}
    \overline{\rho\bu}=\overline{\rho}\bU,
\end{align*}
we apply the same argument to
$\left|\int_{[0,T]\times\O}\phi(\overline{\rho\bu}-\overline{\rho}\bU)\dx\dt\right|$, yielding,
\begin{align*}
    &\left|\int_{[0,T]\times\O}\phi(\overline{\rho\bu}-\overline{\rho}\bU)\dx\dt \right|\leq \int_{[0,T]\times \O}\phi \lan\nu_{t,x},\left| s\beta-s^2\bU\right|\ran \dx\dt+\int_{[0,T]\times \O}\phi \lan\nu_{t,x}^{\infty},\left| s_{\infty}\beta_{\infty}-s^2_{\infty}\bU\right|\ran \dl_t(x)\dt,\\
    &\leq \int_{[0,T]\times \O}\phi \lan\nu_{t,x},c\d s^2+C_\d\left| \beta-s\bU\right|^2\ran \dx\dt+\int_{[0,T]\times \O}\phi \lan\nu_{t,x}^\infty,c\d s^2_\infty+C_\d\left| \beta_{\infty}-s_{\infty}\bU\right|^2\ran \dl_t(x)\dt,\\
    &=c\d\int_{[0,T]\times \O}\phi \overline{\rho} \dx\dt.
\end{align*}
Finally to achieve  
\begin{align*}
    \overline{\rho\bu\otimes\bu}=\overline{\rho}\bU\otimes\bU,
\end{align*}
we estimate $\left|\int_{[0,T]\times\R^d}\phi(\overline{\rho u_iu_j}-\overline{\rho}U_iU_j)\dx\dt \right|$ yielding,
\begin{align*}
    &\left|\int_{[0,T]\times\O}\phi(\overline{\rho u_iu_j}-\overline{\rho}U_iU_j)\dx\dt\right|\\
    &\leq \int_{[0,T]\times\O}\phi\lan \nu_{t,x},\left|\b_i\b_j-s^2U_iU_j\right|\ran\dx\dt+\int_{[0,T]\times\O}\phi\lan \nu_{t,x}^\infty,\left|(\b_{\infty})_i(\b_{\infty})_j-s^2_\infty U_iU_j\right|\ran\dl_t(x)\dt,\\
    &\leq \int_{[0,T]\times \O}\phi\lan\nu_{t,x},c\d s^2+C_\d\left(|\b_i-sU_i|^2+|\b_j-sU_j|^2\right)\ran \dx \dt\\
    & \qquad +\int_{[0,T]\times \O}\phi\lan\nu_{t,x}^\infty,c\d s^2_\infty+C_\d\left(|(\b_\infty)_i-s_\infty U_i|^2+|(\b_\infty)_j-s_\infty U_j|^2\right)\ran \dl_t(x) \dt,\\
    &=c\d \int_{[0,T]\times \O} \phi \overline{\rho} \dx \dt.
\end{align*}
Thus we can conclude \eqref{ident1}, \eqref{ident2}, and \eqref{ident3} all hold almost everywhere in $(0,T)\times \Omega$.

With this in hand we are ready to make the last identification of $\overline{\rho}=r$, almost everywhere in $(0,T)\times \Omega$.

As the pairs $(\overline{\rho},\bU)$ and $(r,\bU)$ both satisfy the transport equation with strong flow $\bU$, we apply Proposition \ref{ineqFK} to get
\begin{align*}
    W_2^2(\overline{\rho},r)\leq C\cE^{mv}_{rel}(\tau)+C\|\overline{\rho}_0-r_0\|_{TV}+Cc(\O)\||x|^2(\overline{\rho}_0-r_0)\|_{TV}=0,
\end{align*}
a.e. in $(0,T) \times \O$.
\end{proof}

\begin{remark}
    Since the density cannot be decoupled from the velocity variable in the relative entropy inequality, the measure family of probabilities $\mu_{t,x}$ corresponding to the variables $(\sqrt{\rho},\sqrt{\rho}\bu)$ satisfies
    \begin{align*}
        \mathrm{supp}(\pi_2\mu)_{t,x}\subset \{\bv: \bv=\a\bU(t,x), \a \in [0,\infty)\},
    \end{align*}
    so that the projection of the generalized Young measure lies on the one-dimensional subspace determined by the direction of the strong velocity $\bU$.
\end{remark}

\section*{Acknowledgements}
\noindent Piotr Gwiazda was supported by National Science Center, Poland through project no. 2023/51/B/ST1/01546. David N. Reynolds and Agnieszka Świerczewska-Gwiazda were supported by National Science Center, Poland through project no. 2023/50/A/ST1/00447. Jakub Woźnicki was partially supported by the Simons Foundation grant (award no. SFI-MPS-T-Institutes-00010825) and from State Treasury funds as part of a task commissioned by the Minister of Science and Higher Education under the project “Organization of the Simons Semesters at the Banach Center - New Energies in 2026-2028” (agreement no. MNiSW/2025/DAP/491).”

\appendix\section{$L^2(I)$ convergence for generating sequences}\label{App1}

Here we provide a necessary result for closing the relative entropy argument when the nonlocal alignment force is present. We construct sequences of functions which generate a space-time generalized Young measure, where the approximating functions need not solve an approximating scheme for the PDE, contrast with the construction in \cite{SzWi2012}. Meanwhile, tracking the density projection throughout the approximation, conservation of mass for the measure-valued solution guarantees an $L^2$ in time convergence in the space of measures under the flat metric. This allows us to retain the key aspect of the weak-strong uniqueness result as this generation result holds for any generalized Young measure and thus their is no dependence on any approximating sequence which is usually used for proving existence of solutions. In particular, letting $D=I\times \Omega$ for $I=[0,T]$ a bounded interval and $\Omega=\T^d$ the periodic spatial domain, we obtain the existence of a sequence of functions which generates the generalized Young measure while also providing strong convergence of the measure-valued density function $\overline{\rho}_t$ in the metric space $L^2(0,T;\cM^+(\O),d_f)$.

First we define the spatial measure
\begin{align*}
    \overline{\rho}_t:=\lan \nu_{t,x},s^2\ran \dx+\lan \nu^\infty_{t,x},s_{\infty}^2 \ran \dl_t(x),
\end{align*}
which is in $L^2(0,T;\cM^+(\T^d),d_f)$ as
\begin{align*}
    \int_I \overline{\rho}_t(\T^d)^2 \dt=\int_I \left(\int_\Omega\lan \nu_{t,x},s^2\ran \dx+\lan \nu^\infty_{t,x},s_{\infty}^2 \ran \dl_t(x)\right)^2 \dt=|I|M^2<\infty,
\end{align*}
by the conservation of mass.

\begin{lemma}\label{Appendix1}
     Let $D=I\times \O$ represent the time-space domain with $\O=\T^d$ the periodic spatial domain. Suppose that the triple $(\nu,\l,\nu^\infty)$ satisfies 
     
     \begin{itemize}
        \item $\nu\in L^\infty_w(D; \cP(\R_+\times \R^d\times \R^d))$,
        \item $\lambda\in \cM^+(\overline{D})$,
        \item $\nu^\infty\in L^\infty_w((D, \lambda); \cP(\S_2^+))$,
    \end{itemize}
    where
    $$
    \S^+_2 := \{(s_\infty, \beta_\infty, F_\infty)\in \R^{2d+1},\, |s_\infty|^4 + |\beta_\infty|^4 + |F_\infty|^4= 1,\, s_\infty\geq 0\},
    $$
     and further assume the disintegration $\lambda(\diff t, \diff x)=\l_t(\diff x) \diff t$ of the concentration measure with respect to Lebesgue measure on $I$. Further suppose that the projection $P(s,\b,F)=s^2$ satisfies
     \begin{align*}
         \int_{\T^d}\lan \nu_{t,x}, s^2 \ran \dx +\int_{\T^d}\lan \nu_{t,x}^\infty,s_\infty^2\ran \dl_t(x)=M, \quad \text{a.e. } t\in I
     \end{align*}
     Then there exists a sequence of functions $\vartheta_n=(s_n,\beta_n,F_n):D\to\R_+\times \R^d\times \R^d$ such that $\vartheta_n$ generates $(\nu, \lambda, \nu^\infty)$ and
    \begin{align}\label{conv:L_2_step_function_generation}
    \mu_t^n=s_n^2\diff x \rightarrow \overline{\rho}_t ,\text{ strongly in }L^2(I; (\mathcal{M}^+(\T^d), d_f)).
    \end{align}
\end{lemma}

\begin{proof}

\textbf{\underline{Step 1: Spatial Averaging.}} For $l\in \N$, consider a partition of the  spatial domain $\T^d$ into a collection of cubes $\{\mathcal{C}^l_i\}_i$ of side length $\frac{1}{l}$ with vertices on the lattice $\frac{1}{l}\Z^{d} \cap \T^d$.  Further suppose the cubes contain their boundary on  the left and top sides and open on the right and bottom sides so that they completely tile the spatial domain. Define
    \begin{align*}
    \langle {\nu}^l_{t,x}, f\rangle &= \frac{1}{|\mathcal{C}^l_i|}\int_{\mathcal{C}^l_i}\langle\nu_{t,y}, f\rangle\diff y,\,x\in \mathcal{C}^l_i, \ t \in I,\\
     {\lambda}^l_t &= \sum_{i}\frac{\lambda_t(\mathcal{C}^l_i)}{|\mathcal{C}^l_i|}\mathcal{L}^{d}|_{\mathcal{C}^l_i}, \ t \in I,\\
    \langle {\nu}^{l, \infty}_{t, x}, g\rangle &= \frac{1}{\lambda_t(\mathcal{C}^l_i)}\int_{\mathcal{C}^l_i}\langle\nu^\infty_{t,y}, g\rangle\diff\lambda_t(y),\, x\in \mathcal{C}^l_i, \ t \in I, \ \l_t(\cC_i^l)>0
    \end{align*}
    so that $ {\lambda}^l(\diff t,\diff x)= {\l}^l_t(\diff x)\diff t$ and for  $f\in C_0(\R_+\times \R^d\times \R^d))$ and $g\in C(\S_2^+)$. When $\l_t(\cC_i^l)=0$, there is no concentration effects in that cell and we can take any arbitrary probability measure on $\S_2^+$. We begin with showing the generation part. Take an admissible $f$ and its recession function $f^\infty$. We need to show that
    $$
    \langle {\nu}^l_{t,x}, f\rangle\diff x\diff t + \langle {\nu}^{l, \infty}_{t, x}, f^\infty\rangle\diff {\lambda}^l(t, x) \rightharpoonup^* \langle\nu_{t,x}, f\rangle\diff x\diff t + \langle\nu^{\infty}_{t, x}, f^\infty\rangle\diff\lambda(t, x)
    $$
    in $\mathcal{M}(D)$. To do so fix $\phi\in C(\overline{D})$, then
    \begin{equation*}
        \begin{split}
            &\int_D\phi\langle {\nu}^l_{t,x}, f\rangle\diff x\diff t + \int_D\phi\langle {\nu}^{l, \infty}_{t, x}, f^\infty\rangle\diff {\lambda}^l(t, x),\\
            &\,= \sum_i\int_I\frac{1}{|\mathcal{C}^l_i|}\int_{\mathcal{C}^l_i}\langle\nu_{t, y}, f\rangle\diff y \int_{\cC_i^l}\phi(t,x)\dx\dt+ \sum_i\int_I\frac{1}{|\cC_i^l|}\int_{\cC_i^l}\langle\nu^\infty_{t, y}, f^\infty\rangle\diff\lambda_t(y)\int_{\cC_i^l}\phi(t, x)\diff x\diff t,\\
            &\, \rightarrow \int_D \phi \lan \nu_{t,x},f\ran \dx\dt + \int_D \phi \lan \nu_{t,x}^{\infty}, f^\infty\ran \dl(t,x),
        \end{split}
    \end{equation*}
    where the convergence can be done due to the continuity of $\phi$. We now show that the convergence of the density in $L^2(I;\cM^+(\T^d),d_f)$, is preserved by this approximation. That is for every $t \in I$, we define the spatial measure
    $$ {\mu}_t^l=\langle  {\nu}_{t,x}^l, s^2\rangle \dx+\langle {\nu}^{l, \infty}_{t, x}, s^2_\infty\rangle\diff {\lambda}_t^l,$$ and we wish to show$$
     {\mu}_t^l \rightarrow \overline{\rho}_t, \quad \text{strongly in} \ L^{2}(I;\cM^+(\T^d),d_f).
    $$
    
    Fix an arbitrary $\psi\in C(\T^d)$, such that $\|\psi\|_\infty, \|\nabla_x\psi\|_\infty\leq 1$. Then,
    \begin{equation*}
        \begin{split}
            \int_{\O}\psi(x) \diff ( \mu_t^l-\overline\rho_t)(x) =\sum_{i}\int_{\mathcal{C}^l_i}\psi(x) -\left( \frac{1}{|\mathcal{C}^l_i|}\int_{\mathcal{C}^l_i}\psi(y)\diff y\right)\diff \overline\rho_t(x),
        \end{split}
    \end{equation*}
    since $\|\nabla \psi\|_\infty\leq 1$
    \begin{align*}
        \left|\psi(x) -\left( \frac{1}{|\mathcal{C}^l_i|}\int_{\mathcal{C}^l_i}\psi(y)\diff y\right)\right|\leq \frac{C}{l},
    \end{align*}
    and
    \begin{equation*}
        \begin{split}
            &\sum_{i}\int_{\mathcal{C}^l_i}\left|\psi(x) - \frac{1}{|\mathcal{C}^l_i|}\int_{\mathcal{C}^l_i}\psi(y)\diff y\right|\diff \overline\rho_t(x) \leq \frac{1}{l}\sum_i\int_{\mathcal{C}^l_i}C\|\nabla_x\psi\|_\infty\diff \overline\rho_t(x)\\
            &\quad\leq \frac{C}{l} M.
        \end{split}
    \end{equation*}
    where we again used the conservation of mass $\overline\rho_t(\T^d)=M$.
    Taking the supremum over $\psi$ we deduce
    \begin{equation*}
        \begin{split}
            \int_I d_f^2( \mu_t^l, \overline\rho_t)\diff t \leq \frac{C^2}{l^2}\int_I \overline\rho_t(\T^d)^2\diff t= \frac{C^2}{l^2}|I|M^2\rightarrow 0,\text{ as }l\to +\infty.
        \end{split}
    \end{equation*}

    \textbf{\underline{Step 2: Time Averaging.}} Now we can take $l$ to be fixed so that our measures $(\nu_t^l, \l_t^l, \nu_t^{l,\infty})$ are piecewise constant in space. We now take $k\in \N$, and further partition the temporal part of the domain $I$ into a collection of intervals $\{\mathcal{I}_j^k\}$ of side length $\frac{|I|}{k}$, again closed on the left and open on the right. Let us introduce some notation. Let

    \begin{align*}
        a_i^l(\t):=\frac{\l_\t(\cC_i^l)}{|\cC_i^l|}
    \end{align*}
    the constant value on each spatial cell $\cC_i^l$ used in the previous step and its time-averaged value
    \begin{align*}
        a_{ij}^k:=\frac{1}{|\cI_j^k|}\int_{\cI_j^k}a_i^l(\t) \diff\tau.
    \end{align*}
    
    Then the space-time averaged generalized Young measure is given by the triple
    \begin{align*}
    \langle\nu^k_{t,x}, f\rangle &= \frac{1}{|\mathcal{I}^k_j|}\int_{\mathcal{I}^k_j}\langle \nu_{\t,x}^{l}, f\rangle\diff \t,\, t\in \mathcal{I}^k_j,\\
    \lambda^k_t&=\sum_i a_{ij}^k \cL^{d}|_{\cC_i^l}(x), \ t\in \cI_j^k\\
    \langle\nu^{k, \infty}_{t, x}, g\rangle &= \frac{1}{|\mathcal{I}^k_j|a_{ij}^k}\int_{\mathcal{I}^k_j} a_i^l(\tau)\langle\nu^{l,\infty}_{\t,x}, g\rangle\diff\tau,\, t\in \mathcal{I}^k_j, \, x\in \cC_i^l, \ a_{ij}^k>0,
    \end{align*}
   for $f\in C_0(\R_+\times \R^d\times \R^d))$ and $g\in C(\S_2^+)$.

   We again begin with showing the generation part. Take an admissible $f$ and its recession function $f^\infty$. We need to show that
    $$
    \langle\nu^k_{t,x}, f\rangle\diff x\diff t + \langle\nu^{k, \infty}_{t, x}, f^\infty\rangle\diff\lambda^k(t, x) \rightharpoonup^* \langle\nu^l_{t,x}, f\rangle\diff x\diff t + \langle \nu^{l,\infty}_{t, x}, f^\infty\rangle\diff\lambda^l(t, x)
    $$
    in $\mathcal{M}(D)$. To do so fix $\phi\in C(\overline{D})$, then
    \begin{equation*}
        \begin{split}
            &\int_D\phi\langle\nu^k_{t,x}, f\rangle\diff x\diff t + \int_D\phi\langle\nu^{k, \infty}_{t, x}, f^\infty\rangle\diff \lambda^k(t, x)\\
            &\,= \sum_j\int_{\cI_j^k}\int_{\T^d}\phi(t,x) \frac{1}{|\cI_j^k|}\int_{\cI_j^k}\langle \nu^l_{\t, x}, f\rangle\diff \t \dx\dt+ \sum_j\int_{\cI_j^k}\int_{\T^d}\phi(t,x)\frac{1}{|\cI_j^k|}\int_{\cI_j^k}a_i^l(\tau)\langle \nu^{l,\infty}_{\t, x}, f^\infty\rangle  \diff \tau \diff x\diff t.
        \end{split}
    \end{equation*}
    In the oscillation part, we exchange the time integrals by Fubini and use the continuity of $\phi$ to obtain
    $$\sum_j\int_{\cI_j^k}\int_{\T^d}\phi(t,x) \frac{1}{|\cI_j^k|}\int_{\cI_j^k}\langle \nu^l_{\t, x}, f\rangle\diff \t \dx\dt \to \int_D \phi(t,x) \langle \nu_{t,x}^l, f\rangle \diff x \diff t.$$
    We proceed similarly in the concentration part,
    \begin{align*}\sum_i\sum_j\int_{\cI_j^k}\int_{\cC_i^l}\phi(t,x)\frac{1}{|\cI_j^k|}\int_{\cI_j^k}a_i^l(\tau)\langle \nu^{l,\infty}_{\t, x}, f^\infty\rangle  \diff \tau \diff x\diff t&= \sum_j\int_{\cI_j^k}\int_{\T^d}\int_{\cI_j^k}\frac{\phi(t,x)}{|\cI_j^k|}\dt \ \langle \nu^{l,\infty}_{\t, x}, f^\infty\rangle \diff \l_\t^l(x) \diff \tau, \\
    &\to \int_D \phi(t,x) \langle \nu_{t,x}^{l,\infty}, f^\infty\rangle \diff \l_t^l(x) \diff t,\end{align*}
    completing the generation part.

     We now show that the convergence in $L^2(I;\cM^+(\T^d),d_f)$, is preserved by this approximation.

 For every $t \in \cI_j^k$, we have the spatial measure
$$\mu_t^k=\langle \nu_{t,x}^k, s^2\rangle \dx+\langle\nu^{k, \infty}_{t, x}, s^2_\infty\rangle\diff\lambda_t^k(x), $$ and we wish to show$$
    \mu_t^k \rightarrow \mu_t^l, \quad \text{strongly in} \ L^{2}(I;\cM^+(\T^d),d_f).
    $$

    Note, both $\mu_t^k$ and $\mu_t^l$ have been spatially averaged so that they are constant on each spatial cell $\cC_i^l$. Meanwhile, $\mu_t^k$ is also constant on time intervals $\cI_j^k$.  Thus we can define for $t\in \cI_j^k$, the mass of the density measures of each cell $\cC_i^l$ with respect to each of the spatial measures.
    \begin{align*}
        m_i^l(t):=\mu_t^l(\cC_i^l)&=\int_{\cC_i^l} \lan \nu_{t,x}^l,s^2\ran \dx + \int_{\cC_i^l} \lan \nu_{t,x}^{l,\infty}, s_\infty^2\ran\diff\l_t^l(x),\\
        &=\int_{\cC_i^l} \lan \nu_{t,x}^l,s^2\ran \dx + \int_{\cC_i^l} a_i^l(t)\lan \nu_{t,x}^{l,\infty}, s_\infty^2\ran\diff x, \nonumber
    \end{align*}
    and similarly
    \begin{align*}
        m_{i,j}^{k,l}(t):=\mu_t^k(\cC_i^l)&=\int_{\cC_i^l} \lan \nu_{t,x}^k,s^2\ran \dx + \int_{\cC_i^l} \lan \nu_{t,x}^{k,\infty}, s_\infty^2\ran\diff\l_t^k(x),\\
        &=\int_{\cC_i^l} \frac{1}{|\cI_j^k|}\int_{\cI_j^k}\lan \nu_{\t,x}^l,s^2\ran \dtau \dx + \int_{\cC_i^l}\frac{1}{|\cI_j^k|}\int_{\cI_j^k} a_i^l(\tau)\lan \nu_{\t,x}^{l,\infty}, s_\infty^2\ran\diff\tau \dx, \nonumber\\
        &= \frac{1}{|\cI_j^k|}\int_{\cI_j^k}\int_{\cC_i^l}\lan \nu_{\t,x}^l,s^2\ran \dx\dtau  + \frac{1}{|\cI_j^k|}\int_{\cI_j^k}\int_{\cC_i^l} a_i^l(\tau)\lan \nu_{\t,x}^{l,\infty}, s_\infty^2\ran \dx\diff\tau, \nonumber\\
        &=\frac{1}{|\cI_j^k|}\int_{\cI_j^k} m_i^l(\tau) \dtau,\nonumber
    \end{align*}
where the order of integration was switched using Tonelli's theorem since both the oscillation and concentration parts of the density variable produce nonnegative measurable integrands. Further, the conservation of mass implies that each $m_i^l,m_{i,j}^{k,l}\leq M$ and hence are elements of $L^2(I)$.
We also note that the spatial measures can then be rewritten as
\begin{align}\label{eq:spatialmeasures}
    \mu_t^l=\sum_i \frac{m_i^l(t)}{|\cC_i^l|} \chi_{\cC_i^l}(x)\dx, \qquad \mu_t^k=\sum_i \frac{m_{i,j}^{k,l}(t)}{|\cC_i^l|}\chi_{\cC_i^l}(x) \dx.
\end{align}
We need the following quick lemmas before proceeding.
\begin{lemma}
    Define the time-average operator over the interval $\cI_j^k$ as $A_kf(t):=\frac{1}{|\cI_{j}^k|}\int_{\cI_j^k} f(\tau)\dtau$ when $t\in \cI_j^k$. Then $A_kf$ is a contraction on $L^2(I)$.
\end{lemma}
\begin{proof}
    \begin{align*}
        \|A_kf\|_{L^2(I)}^2&=\sum_{j}|\cI_j^k|\left|\frac{1}{|\cI_j^k|}\int_{\cI_j^k}f(\tau) \dtau\right|^2,\\
        &\leq \sum_j \int_{\cI_j^k} |f(\tau)|^2 \dtau=\|f\|_{L^2(I)}^2,
    \end{align*}
    by Jensen's inequality.
\end{proof}
\begin{lemma}\label{time-average-contraction}
    For $f\in L^2(I)$ the time-average operator satisfies $A_kf \to f$ strongly in $L^2(I)$ as $k\to \infty$
\end{lemma}
\begin{proof}
   Let $\e>0$. As $C(\overline{I})$ is dense in $L^2(I)$ let $g \in C(\overline{I})$ such that $\|f-g\|_{L^2(I)}<\frac{\e}{3}$. Then
   \begin{align*}
       \|A_kf-f\|_{L^2(I)}&\leq \|A_k(f-g)\|_{L^2(I)}+\|A_kg-g\|_{L^2(I)}+\|g-f\|_{L^2(I)},\\
       &\leq \frac{2\e}{3}+\|A_kg-g\|_{L^2(I)},
   \end{align*}
   since $A_k$ is a contraction on $L^2(I)$. Then for $k>K$ large enough $\|A_kg-g\|_{L^2(I)}<\frac{\e}{3}$ by the continuity of $g$, thus
   \begin{align*}
       \|A_kf-f\|_{L^2(I)}<\e.
   \end{align*}
\end{proof}
     Now taking $m_i^l\in L^2(I)$ we get $\int_I |m_{i,j}^{k,l}-m_i^l|^2 \dt \to 0$ as $k \to \infty$, for each $i$. In order to get control over the flat metric distance between the two spatial measures we recall that $\mu_t^k$ and $\mu_t^l$ are comprised of the functions $m_{i,j}^{k,l}$ and $m_i^l$ respectively via \eqref{eq:spatialmeasures}. Since the measures are constant over the same spatial partitions the total variation of the difference between the two measures over the spatial domain is given by
     \begin{align*}
         |\mu_t^k-\mu_t^l|(\T^d)=\sum_i\int_{\cC_i^l}\frac{|m_{i,j}^{k,l}(t)-m_i^l(t)|}{|\cC_i^l|} \dx=\sum_i |m_{i,j}^{k,l}(t)-m_i^l(t)|,
     \end{align*}
and therefore
     \begin{align*}
        \int_Id_f^2(\mu_t^k,\mu_t^l) \dt\leq  \int_I(|\mu_t^k-\mu_t^l|(\T^d))^2\dt &=\int_I\left(\sum_i |m_{i,j}^{k,l}(t)-m_i^l(t)|\right)^2\dt,\\
        &\leq N_l\sum_i^{N_l}\int_I |m_{i,j}^{k,l}(t)-m_i^l(t)|^2 \dt \to 0,
     \end{align*}
     where $N_l$ is the finite fixed number of spatial cells $\cC_i^l$ and Cauchy-Schwarz and the fact that the total variation distance bounds the flat-metric distance were both used. Finally, the convergence of $\int_I |m_{i,j}^{k,l}(t)-m_i^l(t)|^2 \dt \to 0$ is given by Lemma \ref{time-average-contraction} since $m_{i,j}^{k,l}=A_km_i^l$ for each $t\in \cI_j^k$.

    \textbf{\underline{Step 3: Cut-off.}} By the preceding steps, we may assume that $\nu$ and $\nu^\infty$ are piecewise constant on the space-time cells
$$ \mathcal Q_{j,i}^{k,l}:=\mathcal I_j^k\times\mathcal C_i^l,$$ and that $\lambda$ has a piecewise constant density with respect to $(dt,dx)$, in particular given by $a_{ij}^k$. For each space-time cell, denote the corresponding probability measures by
\begin{align*}
    \nu_{j,i}^k\in\mathcal P(\mathbb R_+\times\mathbb R^d\times\mathbb R^d), \qquad \nu_{j,i}^{k,\infty}\in\mathcal P(\mathbb S_2^+).
\end{align*}
     Let $\Psi_R$ be a standard continuous cut-off function, which is equal to $1$ on $B_R(0)$ (for a ball defined in the state space $\R_+\times \R^d\times\R^d$), $0$ outside $B_{R+1}(0)$ and $0\leq \Psi_R\leq 1$, so that $\Psi_R(s,\beta,F)\to 1$ as $R\to \infty$. Define
    $$
    M^R_{j,i} := \langle\nu_{j,i}^k, 1 - \Psi_R\rangle,
    $$
    and
    $$
    \nu_{j,i}^R := \Psi_R\nu_{j,i}^k + M^R_{j,i}\delta_0.
    $$

    Then $\nu_{j,i}^R$ is a probability measure supported on $B_{R+1}(0)$, since
\begin{align*}
    \lan\nu_{j,i}^R,1\ran=\lan\nu_{j,i}^k,\Psi_R\ran+M_{j,i}^R=1.
\end{align*}
    We begin with showing the generation. As in the previous steps, consider a $\phi\in C(\overline{D})$. Note that we only need to consider the convergence of oscillation part, as the concentration part remains unchanged. Taking admissible $f$, we obtain using Lebesgue's dominated convergence theorem
    \begin{align*}
            \sum_{j,i}\int_D\phi\langle\nu^R_{j,i}, f\rangle\diff x\diff t = \sum_{j,i}\int_D\phi\langle\Psi_R\nu_{j,i}^k, f\rangle + \phi\langle M^R_{j,i}\delta_0, f\rangle\diff x\diff t
            \longrightarrow \int_D\phi\langle\nu_{j,i}^k, f\rangle\diff x\diff t.
    \end{align*}
    We move to showing $L^2$ in time convergence. Denote
    \begin{align*}
        \mu_t^R &:= \sum_{i} \chi_{\cC_i^l}(x)\left(\langle\nu^R_{j,i}, s^2\rangle\diff x + \langle\nu^{k,\infty}_{j,i}, s^2_\infty\rangle\diff\lambda_t^k (x)\right), \ \ t\in \cI_j^k.
    \end{align*}
    The concentration part is unchanged, further we have $\lan \d_0,s^2\ran=0$, so that $$\lan \nu_{j,i}^R,s^2\ran-\lan \nu_{j,i}^k,s^2\ran=\lan (\Psi_R-1)\nu_{j,i}^k,s^2\ran.$$
    
    Thus for an arbitrary $\psi\in C(\T^d)$ such that $\|\psi\|_\infty, \|\nabla_x\psi\|_\infty\leq 1$ we have
    \begin{align*}
         \left|\int_{\T^d}\psi(x)\diff (\mu_t^R-\mu_t^k)(x)\right| =\left|\sum_i\int_{\T^d}\psi(x)\chi_{\cC_i^l}(x)(\langle(\Psi_R - 1)\nu_{j,i}^k, s^2\rangle + \langle M^R_{j,i}\delta_0, s^2\rangle)\diff x\right|.
    \end{align*}
 Taking a supremum over $\psi$ 
 \begin{align*}
    \left|\int_{\T^d}\psi(x)\diff (\mu_t^R-\mu_t^k)(x)\right| &\leq \sum_i\int_{\cC_i^l}\langle (1 - \Psi_R)\nu_{j,i}^k, s^2\rangle\diff x,
 \end{align*}
 and using Lebesgue's dominated convergence theorem, again from the conservation of mass, we get
    \begin{align*}
        \int_Id_f^2(\mu_t^R, \mu_t^k)\diff t \leq \int_I\left(\int_{\T^d}\langle(1 - \Psi_R)\nu_{j,i}^k,s^2\rangle \dx\right)^2 \dt \rightarrow 0.
    \end{align*}

\textbf{\underline{Step 4: Density of discrete measures}}
By the previous step, each $\nu_{j,i}^R$ is supported in the compact set $\overline B_{R+1}(0)$. Since finitely supported probability measures are dense in $\cP(B_{R+1}(0))$ in the weak* topology, see \cite{Bauer+2001}, for each cell \(\mathcal Q_{j,i}^{k,l}\) we can approximate by discrete measures, both for the oscillation and the concentration part since $\S_2^+$ is also compact.
\begin{align*}
    \nu_{j,i}^{m}=\sum_{q=1}^{N_m}\a_{ij,q}^{m}\delta_{z_{ij,q}^{m}}, \qquad \nu_{j,i}^{m,\infty}=\sum_{\k=1}^{M_m}\gamma_{ij,\k}^m\delta_{\zeta_{ij,\k}^m}, \qquad \lambda_{j,i}=\lambda^k(t,x)(Q_{j,i}^{k,l})=\l_t^k(\cC_i^l)=a_{ij}^k \dx \ \text{for} \ t\in \cI_j^k,
\end{align*}
with $z_{ij,q}^{m}=(s_{ij,q}^m,\beta_{ij,q}^m,F_{ij,q}^m)\in \R_+\times \R^d\times \R^d,\zeta_{ij,\k}^{m}=((s_\infty)_{ij,\k}^{m},(\beta_\infty)_{ij,\k}^{m},(F_\infty)_{ij,\k}^{m}) \in \S_2^+$  as well as $$\alpha_{ij,q}^{m},\gamma_{ij,\k}^m\geq 0, \quad \sum_{q=1}^{N_m} \alpha_{ij,q}^{m}=1=\sum_{\k=1}^{M_m}\gamma_{ij,\k}^m$$
so that each $\nu_{j,i}^{m}$ and $\nu_{j,i}^{m,\infty}$ are probability measures. Note that the weak* convergence is at the level of $\nu_{j,i}^{m} \rightharpoonup^* \nu_{j,i}^R$ and $\nu_{j,i}^{m,\infty}\rightharpoonup^*\nu_{j,i}^\infty$, as $N_m,M_m \to \infty$, so that for any admissible $f$ with continuous recession function $f^\infty$, the following convergence holds
$$\lan \nu_{j,i}^{m}, f\ran \to \lan \nu_{j,i}^R,f\ran, \qquad \lan \nu_{j,i}^{m,\infty}, f^\infty\ran \to \lan \nu_{j,i}^{k,\infty},f^\infty\ran,$$ for each $i$. This convergence is stronger than what is needed for the generation argument and thus implies the weak* convergence $\lan \nu_{j,i}^{m}, f\ran \dx\dt+\lan \nu_{j,i}^{m,\infty},f^{\infty}\ran \dl^k(t,x) \rightharpoonup^* \lan \nu_{j,i}^R,f\ran \dx\dt+\lan \nu_{j,i}^\infty,f^\infty\ran \dl^k(t,x)$.\\

Now define the spatial measure
\begin{align*}
    \mu_t^m:&= \sum_{i} \chi_{\cC_i^l}(x)\left(\lan \nu_{j,i}^{m},s^2\ran \dx +\lan \nu_{j,i}^{m,\infty},s^2_\infty\ran \dl_t^k(x)\right), \ \ t\in\cI_j^k,\\
    &=\sum_i \chi_{\cC_i^l}(x)\left(\sum_{q=1}^{N_m}\a_{ij,q}^m(s_{ij,q}^m)^2 \dx+\sum_{\k=1}^{M_m}\g_{ij,\k}^m((s_\infty)_{ij,k}^m)^2\dl_t^k(x)\right).
\end{align*}
Now the $L^2(I;(\cM^+(\T^d),d_f))$ convergence is obtained by considering $\psi \in C(\T^d)$ with $\|\psi\|_\infty,\|\nabla \psi\|_\infty\leq 1$ and bounding
\begin{align*}
    \left|\int_{\T^d} \psi(x) \diff (\mu_t^m-\mu_t^R)(x)\right|&=\left|\sum_i\int_{\cC_i^l} \psi(x) \left(\langle \nu_{j,i}^{m}-\nu_{j,i}^R,s^2\rangle \dx +\lan \nu_{j,i}^{m,\infty}-\nu_{j,i}^{\infty},s^2_\infty\ran \diff \l_t^k(x)\right)\right|,\\
    &\leq \sum_i|\cC_i^l||\langle \nu_{j,i}^{m}-\nu_{j,i}^R,s^2\rangle|+|\l_{t}^k(\cC_i^l)||\lan \nu_{j,i}^{m,\infty}-\nu_{j,i}^{\infty},s^2_\infty\ran|,\\ 
    &\leq C_m:=\max_j\sum_i |\mathcal C_i^l|\left|\lan\nu_{j,i}^m-\nu_{j,i}^R,s^2\ran\right|+a_{ij}^k|\mathcal C_i^l|\left|\lan\nu_{j,i}^{m,\infty}-\nu_{j,i}^{k,\infty},s_\infty^2\ran\right|,
\end{align*}
where the bound is obtained by taking supremum over $\psi$ and using the fact that each measure is now constant over each cube $\cC_i^l$ so that for $t\in \cI_j^k$, $\l_t^k(\cC_i^l)$ is constant. Finally we also have $C_m \to 0$ from the weak* convergence and the fact that there are only finitely many space-time cubes $\cQ_{j,i}^{k,l}$. Hence,
$$\int_I d_f^2(\mu_t^{m},\mu_t^R)\dt\leq |I|C_m^2 \to 0.$$

\textbf{\underline{Step 5: Constructing the generating sequence}}

Now in order to construct the generating sequence so that we retain the $L^2$ in time convergence of the spatial measures we consider the spatial reference cell given by $Y=[0,1)^{d}$. We create a sequence in $n$ of partitions of the domain into the oscillation parts $A_{ij,q}^n$ and concentration parts $B_{ij,\k}^n$ so that
\begin{align*}
    Y=\bigcup_{q=1}^{N_m} A_{ij,q}^n \cup \bigcup_{\k=1}^{M_m} B_{ij,\k}^n=A_{ij}^n \cup B_{ij}^n,
\end{align*}
for disjoint measurable sets $A_{ij,q}^n, B_{ij,\k}^n \subset Y$. Then within the particular space-time cell $\cQ_{j,i}^{k,l}$ we have $$a_{ij}^k=\frac{\l_t^k(\cC_i^l)}{|\cC_i^l|}$$ independent of $t\in \cI_j^k$. Thus we prescribe the volumes
\begin{align*}
    |A_{ij,q}^n|=\alpha_{ij,q}^{m}(1-|B_{ij}^n|), \qquad |B_{ij,\k}^n|=\frac{a_{ij}^k}{b_n^4}\g_{ij,\k}^m,
\end{align*}
where $b_n\to \infty$ and $|B_{ij,\k}^n|\to 0$ in such a way that $|B_{ij,\k}^n|b_n^4\to a_{ij}^k\g_{ij,\k}^m$. Now we define the functions on the oscillation and concentration sets
\[\o_{j,i}^n(t,x)=\begin{cases}
    z_{ij,q}^{m}, & \text{if} \ x\in A_{ij,q}^n,\\
    b_n^2\zeta_{ij,\k}^m, & \text{if} \ x \in B_{ij,\k}^n.
\end{cases}\]

We then define the periodic extension $\tilde \o_{j,i}^n(t,y+\xi)=\o(t,y)$ for each $\xi \in \Z^d$.

Now as the spatial partition $\cC_i^l$ has side lengths $\frac{1}{l}$ we let $\e_n=\frac{1}{ln}$. Then tile $\cC_i^l$ with $\e_n(Y+\xi)$, for $\xi\in \Z^{d}$ such that $\e_n\xi\in \frac{1}{l}\T^d$. Now define for $(t,x) \in \cQ_{j,i}^{k,l}$
\begin{align*}
    \vartheta_n^{j,i}(t,x)=\tilde\o_{j,i}^n\left(t,\frac{x}{\e_n}\right).
\end{align*}
Then the generating sequence is given by
\begin{align*}
    \vartheta_n(t,x)=\sum_{i,j} \chi_{\cQ_{j,i}^{k,l}}(t,x)\vartheta_n^{j,i}(t,x).
\end{align*}
We now show generation. Let $\phi \in C(\overline{D})$ and consider an admissible $f$, with continuous recession function $f^{\infty}$, then
\begin{align*}
    \int_D \phi(t,x)f(\vartheta_n) \dx \dt&= \sum_{j,i}\int_{\cI_j^k} \int_{\cC_i^l} \phi(t,x)f(\vartheta_n^{j,i}) \dx\dt,\\
    &=\sum_{j,i} \int_{\cI_j^k} \int_{\cC_i^l}\phi(t,x)f\left(\tilde\omega_{j,i}^n\left(t,\frac{x}{\varepsilon_n}\right)\right) \dx \dt,
\end{align*}
Now the change of spatial variables $y= \frac{x}{\e_n}$ gives
\begin{align*}
\sum_{j,i} &\int_{\cI_{j}^k} \int_{\cC_i^l}\phi(t,x)f\left(\tilde\omega_{j,i}^n\left(t,\frac{x}{\varepsilon_n}\right)\right) \dx \dt=\sum_{j,i}\int_{\cI_j^k}\sum_{\xi}\e_n^{d}\int_{Y}\phi(t,\e_n(\xi+y))f\left(\omega_{j,i}^n\left(t,y\right)\right) \dy \dt.
\end{align*}
Now adding and subtracting $\phi(t,\e_n\xi)$ we control the first term by the continuity of $\phi$ so that
$$|\phi(t,\e_n(\xi+y))-\phi(t,\e_n\xi)|\leq C(\e_n),$$
since $y\in [0,1)^d$. Thus we get 
\begin{align*}
\sum_{j,i}\int_{\cI_j^k}\sum_{\xi}\e_n^{d}\int_{Y}\phi(t,\e_n(\xi+y))f\left(\omega_{j,i}^n\left(t,y\right)\right) \dy \dt& =\sum_{j,i}\int_{\cI_j^k}\sum_{\xi}\e_n^{d}\phi(t,\e_n\xi)\int_{Y}f\left(\omega_{j,i}^n\left(t,y\right)\right) \dy \dt +\e_2(n),
\end{align*}
where $\e_2(n) \to 0$ as $n \to \infty$ by the continuity of $\phi$. Now noticing that since $\xi$ were chosen as the lattice points within $\cC_i^l$, so that $\e_n=\frac{1}{ln}$ perfectly tiles $\cC_i^l$, that $\sum_\xi \e_n^d \phi(t,\e_n\xi)$ is a Riemann sum approximation of $\int_{\cC_i^l}\phi(t,x) \dx$. Thus adding and subtracting $\int_{\cC_i^l}\phi(t,x) \dx$ we get
\begin{align*}
\sum_{j,i}\int_{\cI_j^k}\sum_{\xi}\e_n^{d}\phi(t,\e_n\xi)\int_{Y}f\left(\omega_{j,i}^n\left(t,y\right)\right) \dy \dt&=\sum_{j,i}\int_{\cI_j^k}\int_{\cC_i^l}\phi(t,x)\dx\int_{Y} f(\o_{j,i}^n(t,y)) \dy\dt +\e_3(n),
\end{align*}
where $\e_3(n) \to 0$ again by the continuity of $\phi$ and the convergence of the Riemann sum.

Finally,
\begin{align*}
   \sum_{j,i}&\int_{\cI_j^k}\int_{\cC_i^l}\phi(t,x)\dx\int_Y f(\o_{j,i}^n(t,y)) \dy\dt =\sum_{j,i}\int_{\cI_j^k}\int_{\cC_i^l} \phi(t,x)\left(\sum_{k=1}^{N_m}\int_{A_{ij,q}^n} f(z_{ij,q}^m) \dy+\sum_{\k=1}^{M_m} \int_{B_{ij,\k}^n} f(b_n^2\zeta_{ij,\k}^m) \dy\right) \dx\dt,\\
   &=\sum_{j,i}\int_{\cI_j^k}\int_{\cC_i^l} \phi(t,x)\left(\sum_{q=1}^{N_m}|A_{ij,q}^n| f(z_{ij,q}^m) +\sum_{\k=1}^{M_m} |B_{ij,\k}^n| f(b_n^2\zeta_{ij,\k}^m) \right) \dx\dt,\\
   &=\sum_{j,i}\int_{\cI_j^k}\int_{\cC_i^l} \phi(t,x)\left(\sum_{q=1}^{N_m}\a_{ij,q}^m(1-|B_{ij}^n|)) f(z_{ij,q}^m) +\sum_{\k=1}^{M_m} \frac{a_{ij}^k}{b_n^4} \g_{ij,\k}^mf(b_n^2\zeta_{ij,\k}^m) \right) \dx\dt,\\
   &\to \sum_{j,i}\int_{\cI_j^k}\int_{\cC_i^l} \phi(t,x)\left(\sum_{q=1}^{N_m}\alpha_{ij,q}^m f(z_{ij,q}^m) +a_{ij}^k\sum_{\k=1}^{M_m} \g_{ij,\k}^m f^{\infty}(\zeta_{ij,\k}^m) \right) \dx\dt,\\
   &=\int_D \phi(t,x) \sum_{j,i}\left(\lan \nu_{j,i}^{m}, f\ran \dx\dt+\lan \nu_{j,i}^{m,\infty}, f^\infty\ran \dl^k(t,x)\right),
\end{align*}
where the convergence step saw $|B_{ij}^n|\to 0$ and $\frac{f(b_n^2\zeta_{ij,\k}^m)}{b_n^4}\to f^{\infty}(\zeta_{ij,\k}^m)$ as $n\to \infty$.

Now, we need to prove the $L^2$ convergence for the spatial measure generated by $\vartheta_n$. Defining the density projection $$P:\R_+\times \R^d\times \R^d\to \R^+, \quad P(s,\b,F)=s^2,$$
we get the approximating density functions given for $(t,x) \in \cQ_{j,i}^{k,l}$,
\begin{align*}
\mu_t^n(x)&=P(\vartheta_n(t,x))\dx.
\end{align*}
Further, the target density is
\begin{align*}
    \mu_t^m=c_{ij}^m \dx,
\end{align*}
where $$c_{ij}^m:=\sum_{q=1}^{N_m}\a_{ij,q}^m(s_{ij,q}^m)^2 +a_{ij}^k\sum_{\k=1}^{M_m}\g_{ij,\k}^m((s_\infty)_{ij,k}^m)^2.$$

Now let $Q_\xi^n=\e_n(Y+\xi)$ be the microscopic cells contained in a spatial cell $\cC_i^l$. Then with the change of variables $x=\e_n(\xi+y)$, since $Q_\xi^n$ perfectly tiles $\cC_i^l$,
\begin{align*}
    \mu_t^n(Q_\xi^n)&=\int_{Q_\xi^n} P(\vartheta_n(t,x))\dx,\\
    &=\e_n^d\int_YP(\o_{j,i}^n(t,y)\dy,\\
    &=\e_n^d\left(\sum_{q=1}^{N_m}\a_{ij,q}^m(1-|B_{ij}^n|)) (s_{ij,q}^m)^2 +a_{ij}^k\sum_{\k=1}^{M_m}  \g_{ij,\k}^m ((s_\infty)_{ij,\k}^m)^2 \right),\\
    &=\e_n^dc_{ij}^m-\e_n^d|B_{ij}^n|\sum_{q=1}^{N_m}\a_{ij,q}^m(s_{ij,q}^m)^2.
\end{align*}

Meanwhile $$\mu_t^m(Q_\xi^n)=\e_n^dc_{ij}^m,$$ and therefore
\begin{align*}
    |\mu_t^n(Q_\xi^n)-\mu_t^m(Q_\xi^n)|=\e_n^d||B_{ij}^n|\sum_{q=1}^{N_m}\a_{ij,q}^m(s_{ij,q}^m)^2\to 0. 
\end{align*}
Then letting $\psi\in C(\T^d)$ such that $\|\psi\|_\infty,\|\n\psi\|_\infty\leq 1$ we consider for some $x_\xi^n\in Q_\xi^n$,
\begin{align*}
   \left| \int_{Q_\xi^n}\psi(x)\diff (\mu_t^n-\mu_t^m)(x)\right|&\leq\left|\int_{Q_{\xi}^n}(\psi(x)-\psi(x_\xi^n))\diff (\mu_t^n-\mu_t^m)(x)\right|+\left|\psi(x_\xi^n)(\mu_t^n(Q_\xi^n)-\mu_t^m(Q_\xi^n))\right|,\\
    &\leq C\e_n(\mu_t^n(Q_\xi^n)+\mu_t^m(Q_\xi^n))+\e_n^d||B_{ij}^n|\sum_{q=1}^{N_m}\a_{ij,q}^m(s_{ij,q}^m)^2.
\end{align*} 
Summing over $\xi$ gives 
\begin{align*}
   \left| \int_{\cC_i^l} \psi(x) \diff (\mu_t^n-\mu_t^m)\right|\leq C\e_n(\mu_t^n(\cC_i^l)+\mu_t^m(\cC_i^l))+|\cC_i^l||B_{ij}^n|\sum_{q=1}^{N_m}\a_{ij,q}^m(s_{ij,q}^m)^2,
\end{align*}
and finally summing over each $i$ gives the bound for the whole space
\begin{align*}
    d_f(\mu_t^n,\mu_t^m)\leq C\e_n(\mu_t^n(\T^d)+\mu_t^m(\T^d))+\sum_i|\cC_i^l||B_{ij}^n|\sum_{q=1}^{N_m}\a_{ij,q}^m(s_{ij,q}^m)^2\to 0,
\end{align*}
as $n\to \infty$ since there are only finitely many cells and $|B_{ij}^n|\to 0$.

Finally,
\begin{align*}
    \int_I d_f^2(\mu_t^n,\mu_t^m) \dt\to 0,
\end{align*}
by the conservation of mass and Lebesgue's dominated convergence.

Finally a diagonal argument yields the sequence of functions $\vartheta_n$ which generates the original generalized Young measure while also retaining the $L^2(I)$ convergence of the spatial measures corresponding to the density projection.

\end{proof}

\bibliographystyle{plain}   
\bibliography{references}

\end{document}